\documentclass{amsart}

\usepackage[utf8]{inputenc}
\usepackage{amsthm,amsmath,amssymb,latexsym,stmaryrd,lineno}
\usepackage{mathrsfs}
\usepackage{hyperref}
\hypersetup{urlcolor=blue, citecolor=red, linkcolor=blue}

\allowdisplaybreaks
\numberwithin{equation}{section}
\theoremstyle{plain} 

\newtheorem{thm}[equation]{Theorem}
\newtheorem*{thm*}{Theorem}
\newtheorem{prop}[equation]{Proposition}
\newtheorem{coro}[equation]{Corollary}
\newtheorem{lem}[equation]{Lemma}
\theoremstyle{definition} 
\newtheorem{ex}[equation]{Example}

\newtheorem{rem}[equation]{Remark}

\newtheorem{defn}[equation]{Definition}

\renewcommand{\to}{\longrightarrow}
\renewcommand{\mapsto}{\longmapsto}
\newcommand{\Id}{ \mathrm{Id}}
\newcommand{\Mat}{\mathrm{Mat}}

\numberwithin{equation}{section}

\begin{document}

\title{On the decomposability of bilinear spaces of dimension four}

\author{Grégory Berhuy} 

\address{Université Grenoble Alpes\\
  Institut Fourier\\
  CS 40700, 38058 Grenoble cedex 9} 
\email{gregory.berhuy@univ-grenoble-alpes.fr}

\begin{abstract}
In this paper, we study the problem of decomposability of bilinear spaces of dimension four without symmetry, as well as the problem of decomposability of split central simple algebras of degree four with an anti-automorphism. In particular, we show that, contrary to the case of symmetric or skew-symmetric bilinear spaces, these two problems are not equivalent. We will also prove that cohomological invariants do not detect decomposability of bilinear spaces of dimension four
in general, whereas the determinant does for split central simple algebras of degree four with an anti-automorphism.
\end{abstract}

\maketitle

%\linenumbers


%\MSC[2020] 11E72 \sep  11E39

\tableofcontents

\section{Introduction}

{\bf Convention. }In this paper, $F$ is a field, $\overline{F}$ is a fixed algebraic closure of $F$, and $F$-algebras and $F$-vector spaces are implicitely supposed to be finite-dimensional over $F$, unless specified otherwise.

Let $F$ be a field of characteristic different from two.

Central simple $F$-algebras with involution are important objects to study, since they play a crucial 
role in the classification of linear algebraic groups over an arbitrary field. 
Moreover, they may be seen as a natural generalization 
of symmetric/alternating/hermitian spaces (see \cite{KMRT}).

We will say that a central simple $F$-algebra with involution $(A,\sigma)$ is \textit{decomposable} if there exist two central simple $F$-algebras with involutions $(A_1,\sigma_1)$ and $(A_2,\sigma_2)$, with $\deg(A_i)\geq 2, i=1,2$, such that \[(A,\sigma)\simeq (A_1,\sigma_1)\otimes_F (A_2,\sigma_2).\]

Even if the classification of central simple algebras with involution is out of reach, various decomposability criteria, in terms of the vanishing of suitable cohomological invariants, have been obtained by several authors.

For example, decomposability criteria have been proven by Knus, Parimala and Sridharan  for algebras of degree four with an $F$-linear involution in \cite{KPS}, by Karpenko and Qu\'{e}guiner-Mathieu  for algebras of degree four with a unitary involution in \cite{KaQu}, by Qu\'{e}guiner-Mathieu and Tignol for algebras of degree $8$ and $12$ with an orthogonal involutions in \cite{QuTi8} and \cite{QuTi12} respectively, and by Garibaldi, Parimala and Tignol for algebras of degree $8$ with a symplectic involution in \cite{GPT}. More results on decomposability and cohomological invariants may be found in Tignol's survey paper \cite{Ti}.

In this paper, we will focus on the decomposability of central simple algebras of degree four with an $F$-linear anti-automorphism. To motivate things further, let us state the main result of \cite{KPS} precisely.

If  $(A,\sigma)$ is a central simple $F$-algebra with an $F$-linear involution, recall that the {\it determinant} of $(A,\sigma)$ is defined as the square class of any invertible skew-symmetric element.
In \cite{KPS}, Knus, Parimala and Sridharan proved the following result.

\begin{thm*}[\cite{KPS}]
Let $(A,\sigma)$ be a central simple $F$-algebra of degree four with an $F$-linear involution. Then, $(A,\sigma)$ is decomposable if and only if the determinant of $\sigma$ is trivial.
\end{thm*}

Note that this result may be viewed as a perfect analogue of the  much more elementary corresponding result for bilinear spaces.  

Say that $(V,b)$ is \textit{decomposable} if there exist two bilinear spaces $(V_1,\varphi_1)$ and $(V_2,\varphi_2)$, with $\dim_F(V_i)\geq 2$, such that \[ (V,b)\simeq (V_1,\varphi_1)\otimes_F (V_2,\varphi_2). \]

It is known that if $(V,b)$ is a non-degenerate symmetric or skew-symmetric $F$-bilinear space of dimension four, then $(V,b)$ is decomposable if and only if $\det(b)$ is trivial.

In fact, this last result is exactly the result of Knus, Parimala and Sridharan applied to the case $A=\mathscr{L}(V)$, the $F$-algebra of endomorphisms of $V$.
To explain why, let us introduce some notation.

Let $(V,b)$ be a non-degenerate $F$-bilinear space. If $f\in\mathscr{L}(V)$ is an endomorphism of $V$, there is a unique endomorphism $\sigma_b(f)$ such that \[ b(f(x),y)=b(x,\sigma_b(f)(y)) \ \mbox{ for all }x,y\in V . \]

The map $\sigma_b$ is an anti-automorphism of the $F$-algebra $\mathscr{L}(V)$, which is said to be {\it adjoint to $b$}. When $b$ is symmetric or skew-symmetric, this is even an $F$-linear involution. Conversely, any $F$-linear anti-automorphism of $\mathscr{L}(V)$ is adjoint to some non-degenerate bilinear form $b:V\times V\to F$, which is uniquely determined up to multiplication by an element of $F^\times$.

We then have the following result, which is due to Shapiro (\cite[Chapter 9]{Sha}), and which has been generalized by Becher in \cite{Bec}.
 
\begin{thm*}
Let $(V,b)$ be a symmetric or skew-symmetric non-degenerate $F$-bilinear space of dimension four. The following properties are equivalent : 

\begin{enumerate}
    \item the bilinear space $(V,b)$ is decomposable

    \item the central simple $F$-algebra with involution $(\mathscr{L}(V),\sigma_b)$ is decomposable.
\end{enumerate}
\end{thm*}

Note that, while the implication $(1)\Longrightarrow (2)$ is obvious, the reverse implication is not trivial. If $(\mathscr{L}(V),\sigma_b)\simeq (Q_1,\sigma_1)\otimes_F (Q_2,\sigma_2)$ for some quaternion $F$-algebras $Q_1,Q_2$,  there is no reason why we may choose $Q_1$ and $Q_2$ to be split.

In \cite{CT}, Cortella and Tignol extended the notion of determinant to the case of central simple $F$-algebras of even degree with an $F$-linear anti-automorphism $(A,\sigma)$. When $(A,\sigma)=(\mathscr{L}(V),\sigma_b)$, the determinant of $\sigma$ coincides with $\det(b)$ (see \cite[Proposition 9]{CT}). We will not define precisely the determinant in the general case, but let us mention that if $(A,\sigma)$ is decomposable of degree four, then $\det(\sigma)$ is trivial. Indeed, the arguments presented in the proof of Proposition \ref{adjclosed} may be immediately adapted to this case.

Therefore, one could ask whether the result of Knus, Parimala and Sridharan extends to the case of $F$-linear anti-automorphisms. At this point, the answer is no, because of the existence of an additional obstruction. Indeed, since $F$-linear anti-automorphisms over split central simple $F$-algebras are adjoint to bilinear forms, if $(A,\sigma)$ is decomposable, then it is adjoint over $\overline{F}$ to a decomposable $\overline{F}$-bilinear space. In fact, one can say even more. We will show that if $(Q,\tau)$ is an $F$-quaternion algebra with an $F$-linear anti-automorphism, then $(Q,\tau)_{\overline{F}}$ is adjoint to an $F$-bilinear form which is defined over $F$ (see Lemma \ref{quat}).
In particular, if $(A,\sigma)$ is decomposable (where $A$ has degree four), then $(A,\sigma)_{\overline{F}}\simeq (\mathscr{L}(V),\sigma_b)_{\overline{F}}$ for some decomposable $F$-bilinear space $(V,b)$. Note that this condition is automatically fulfilled if $\sigma$ is an arbitrary (not necessarily decomposable) $F$-linear involution. Indeed, in this case, $\sigma$ is adjoint over $\overline{F}$ to a symmetric/skew-symmetric hyperbolic space of dimension four, which is decomposable and defined over $F$. More precisely,  
$(A,\sigma)$ and $(\mathscr{L}(F^2),\sigma_{\langle 1,1\rangle})\otimes_F (\mathscr{L}(F^2),\sigma_h)$ are isomorphic over $\overline{F}$, where $h:F^2\times F^2\to F^2$ is the standard hyperbolic symmetric/skew-symmetric form. 

To sum up, the right problem to consider is in fact the following one.

{\bf Problem. }Let $(A,\sigma)$ be a central simple $F$-algebra of degree four with an $F$-linear anti-automorphism. Assume that  $(A,\sigma)_{\overline{F}}\simeq (\mathscr{L}(V),\sigma_b)_{\overline{F}}$ for some decomposable $F$-bilinear space $(V,b)$. Is it true that $(A,\sigma)$ is decomposable if and only $\sigma$ has trivial determinant ?

In this paper, we will investigate this problem in the split case, which will already happen to be non-trivial.
In this particular case, the problem may be reformulated in a nicer way.

To explain how, we need to recall the notion of asymmetry of a bilinear space.
Let $(V,b)$ be a  non-degenerate $F$-bilinear space. Then, 
there is a unique endomorphism $a_b$ of $V$, called the \textit{asymmetry} of $b$, such that \[ b(x,y)=b(y, a_b(x)) \ \mbox{ for all }x,y\in V. \]
It is easy to see that the asymmetry of $(V_1,\varphi_1)\otimes_F(V_2,\varphi_2)$ is $a_{\varphi_1}\otimes a_{\varphi_2}$, and that
isomorphic/similar bilinear spaces have conjugate asymmetries.
Moreover, one may show that two non-degenerate $F$-bilinear spaces are similar over $\overline{F}$ if and only if they have conjugate asymmetries over $F$ (Lemma \ref{symfbar}).

Let us say that a bilinear space $(V',b')$ \textit{has a decomposable asymmetry} if there exists a decomposable bilinear space $(V,b)$ such that $a_{b'}$ is conjugate to $a_b$.  It follows from the previous considerations that, given a non-degenerate $F$-bilinear space $(V',b')$, there exists a decomposable $F$-bilinear space $(V,b)$ such that $(\mathscr{L}(V'),\sigma_{b'})_{\overline{F}}\simeq (\mathscr{L}(V),\sigma_{b})_{\overline{F}}$   
if and only if $(V',b')$ has a decomposable asymmetry. 
Taking into account that $\det(\sigma_{b'})=\det(b')$, our problem in the split case then may be rephrased as follows.

\textbf{Question 1. }Let $(V',b')$ a non-degenerate $F$-bilinear space of dimension four. Assume that $(V',b')$ has a decomposable asymmetry.

Is is true that $(\mathscr{L}(V'),\sigma_{b'})$ is decomposable if and only if $b'$ has trivial determinant ?

A natural strategy to tackle this problem is to work at the level of bilinear forms. Note that a decomposable bilinear space has a decomposable asymmetry and a trivial determinant. Thus, one may ask the following stronger question.

\textbf{Question 2. }Let $(V',b')$ a non-degenerate $F$-bilinear space of dimension four with a decomposable asymmetry. Is is true that $(V',b')$ is decomposable if and only if $\det(b')$ is trivial ?

Obviously, a positive answer to Question 2 yields a positive answer to Question 1. This approach also leads to the following related problem, which is of independent interest.

\textbf{Question 3. }Let $(V',b')$ be a non-degenerate $F$-bilinear space of dimension four. Is $(V',b')$ decomposable
if and only if $(\mathscr{L}(V'),\sigma_{b'})$ is decomposable ?

In this paper, we will show that Question 1 has a positive answer (Theorem \ref{thmadjoint}). However, if $(V',b')$ has a generic decomposable asymmetry (see Definition \ref{genasym}), we will see that Question 2 has a negative answer in general (Theorem \ref{thmnondec}). In particular, Question 3, quite surprisingly, has a negative answer as well.

This could be the end of the story. However, in order to have a complete picture of the situation, it is very tempting to ask if there are other invariants of bilinear spaces which might detect decomposability in dimension four.

To make this question more precise, let us indicate that, if $\mathbf{GO}(b)$ is the algebraic group of similitudes of $(V,b)$, one may show that $H^1(F, \textbf{GO}(b))$ is in one-to-one correspondence with the set of similarity classes of bilinear spaces $(V',b')$ such that $a_b$ and $a_{b'}$ are conjugate (Proposition \ref{cohogo}). The determinant then may be reinterpreted as a normalized cohomological invariant of degree $1$  with coefficients in $\mathbb{Z}/2\mathbb{Z}$ of the algebraic group $\textbf{GO}(b)$. Despite the negative answer to Question 2, one may wonder if there exists a family  of normalized cohomological invariants of higher degrees which may detect decomposability. 
Note that decomposability is a property which only depends on the similarity class of a bilinear space. Hence, one may ask :

\textbf{Question 4. } Let $(V,b)$ be a decomposable $F$-bilinear space of dimension four. 

Does there exist a family $\mathscr{F}$ of cohomological invariants of $\textbf{GO}(b)$ with coefficients in $\mathbb{Z}/2\mathbb{Z}$ such that, for every field extension $E/F$, and for every $E$-bilinear space $(V',b')$ whose asymmetry is conjugate to $a_b$ over $E$, we have $\eta_E(b')=0$ for all $\eta\in\mathscr{F}$ if and only if $b'$ is decomposable ? 

Alas, in the case of bilinear spaces with a generic asymmetry, we will provide examples of fields $F$ for which the answer is no (Remark \ref{negquest2}, Example \ref{laurentr} and Example \ref{laurent}). However, we will construct an invariant of $\mathbf{GO}(b)$ with values in an appropriate quotient of the Brauer group which detects decomposability in this case (Proposition \ref{dec-detect}).

We now describe the organization of this paper. In Section \ref{sec-lemmas}, we provide several lemmas on asymmetries which will be useful in the sequel. In Section \ref{sec-corres}, we establish an explicit correspondence between similarity classes of bilinear spaces with a fixed asymmetry (up to conjugacy) and the set of equivalence classes of invertible symmetric elements of some algebra with involution with respect to a certain equivalence relation.
We then answer Questions 1 and 2 in the case of bilinear spaces with a non-generic decomposable asymmetry. The  case of a  generic decomposable asymmetry, which is the hardest case, is handled in Section \ref{sec-det-gen}.
In Section \ref{sec-coh}, we reinterpret our object of study in terms of Galois cohomology. We then address Question 4 in Section \ref{sec-cohinv}.

{\bf Notation. }Let $V$ be a finite-dimensional $F$-vector space,  and let $\mathscr{B}$ be an $F$-basis of $V$. If 
$f$ is an endomorphism of $V$,  its matrix representation  with respect to $\mathscr{B}$ will be denoted by $\mathrm{Mat}(f;\mathscr{B})$. Similarly, if $b$ is an $F$-bilinear form on $V$, its matrix representation  with respect to $\mathscr{B}$ will be denoted by $\mathrm{Mat}(b;\mathscr{B})$.

\section{Some lemmas on asymmetries}\label{sec-lemmas}

In this section, we collect some preliminary results on asymmetries which will be useful in the sequel.

\begin{lem}\label{abinv}
Let $F$ be a field of arbitrary characteristic, and let $(V,b)$ be a non-degenerate $F$-bilinear space of dimension $n\geq~1$.
Then $a_b$ is conjugate to its own inverse. In particular, we have \[\mu_{a_b}=\mu_{a_b}(0)^{-1}X^{\deg(\mu_{a_b})}\mu_{a_b}(X^{-1}) \ \mbox{ and } \chi_{a_b}=\chi_{a_b}(0)^{-1}X^n\chi_{a_b}(X^{-1}).\]
\end{lem}

\begin{proof}
Let $\mathscr{B}$ be a fixed basis of $V$, and let $B=\Mat(b;\mathscr{B})$.
Then $\Mat(a_b;\mathscr{B})=B^{-1}B^t$. 
Now, \[  B(B^{-1}B^t)B^{-1}=B^tB^{-1}=(B^{-1}B^t)^{-t}. \]
But, it is known that any matrix is conjugate to its transpose, so $B^{-1}B^t$ is conjugate to its inverse. The rest of the lemma comes from the fact that conjugate endomorphisms have same minimal and characteristic polynomials.
\end{proof}

\textbf{Notation. }Let $\mathscr{E}$ be the canonical basis of $F^2$. If $a\in F\setminus\{2\}$, we denote by $b_a$ the unique $F$-bilinear form on $F^2$ such that \[ \mathrm{Mat}(b_a;\mathscr{E})=\begin{pmatrix}
       1 & 1 \cr 0 & \frac{1}{2-a}
   \end{pmatrix}. \]  
   
 We then have   \[ \mathrm{Mat}(a_{b_a}; \mathscr{E})= A_a\overset{\mathrm{def}}{=}\begin{pmatrix}
 a-1 & -1 \cr 2-a & \phantom{-}1
 \end{pmatrix}. \]

We also denote by $h$ the unique alternating $F$-bilinear form on $F^2$ such that \[ \mathrm{Mat}(h; \mathscr{E})=H\overset{\mathrm{def}}{=}\begin{pmatrix}
    \phantom{-}0 & 1 \cr -1 & 0 
\end{pmatrix}. \]
   
 \begin{lem}\label{lemba}
  Let $(V,b)$ be a non-degenerate bilinear space of dimension two over a field $F$ of arbitrary characteristic. Then 

  \begin{enumerate}
  
  \item $b$ is symmetric if and only if $\mu_{a_b}=X-1$.

  \item $b$ is skew-symmetric if and only if $\mu_{a_b}=X+1$. 

   \item Assume that $b$ is not symmetric nor alternating. 
   
   Then $\mu_{a_b}=\chi_{a_b}=X^2-aX+1$, for some $a\neq 2$. Moreover, $b$ is similar to $b_a$.
  \end{enumerate}
 \end{lem}
\begin{proof}
 Items  $(1)$ and $(2)$ are clear.

Assume now that $b$ is not symmetric nor alternating. Then $\mu_{a_b}$ has degree two.
Otherwise, the last part of Lemma \ref{abinv} implies that $\mu_{a_b}=X\pm 1$. If $F$ has characteristic different from two, the two first items imply that $b$ is symmetric or alternating (since a bilinear form is then skew-symmetric if and only if it is alternating). If $F$ has characteristic two, we have $\mu_{a_b}=X-1$, so $a_b=\Id_V$, meaning that $b$ is symmetric. 

Since $\mu_{a_b}$ has degree two, we get $\chi_{a_b}=\mu_{a_b}=X^2-aX+1$, for some $a\in F$, again by applying the last part of Lemma \ref{abinv}.

We now prove that $b$ is similar to $b_a$, and that $a\neq 2$. Since $b$ is non-alternating, there exists $e_1\in V$ such that $b(e_1,e_1)=\alpha\neq 0$. 

Now, the linear form $b(_-,e_1)$ is non-zero, so its kernel has dimension one. Let $e_2\in V\setminus\{0\}$ such that $b(e_2,e_1)=0$.
Note that $(e_1,e_2)$ is a basis of $V$ ($e_2$ cannot be proportional to $e_1$ since $b(e_2,e_1)=0$ while $b(e_1,e_1)\neq 0$).

We have $\beta=b(e_1,e_2)\neq 0$, since otherwise the matrix representation of $b$ with respect to the basis $(e_1,e_2)$ would be symmetric. Replacing $e_2$
by $\alpha\beta^{-1}e_2$, one may assume further that $b(e_1,e_2)=\alpha$. Setting $c=\alpha^{-1}b(e_2,e_2)$,
we see that the matrix representation of $b$ with respect to $(e_1,e_2)$ is \[ B=\alpha \begin{pmatrix}
  1 & 1 \cr 0 & c  
\end{pmatrix}, \]
with $c\neq 0$ since $b$ is non-degenerate. 
Now $B^{-1}B^t=\begin{pmatrix}
 1-c^{-1} & -1 \cr c^{-1} & \phantom{-}1   
\end{pmatrix}$, and thus $\chi_{a_b}=X^2-(2-c^{-1})X+1$.
Then $a=2-c^{-1}$, and in particular $a\neq 2$.
Now, $c=\frac{1}{2-a}$, and the previous computations show that $b$ is similar to $b_a$.
\end{proof}

\begin{coro}\label{coroba}
Let $(V,b)$ be a non-degenerate bilinear space of dimension two.

Write $\chi_{a_b}=X^2-aX+1$. Then:

\begin{enumerate}
    \item if $a=2$, $b$ is symmetric. In particular, $b$ is similar to $\langle 1,d\rangle$, where $\det(b)=d\in F^\times/F^{\times 2}$ if $b$ is non-alternating, and $b$ is isomorphic to $h$ if $b$ is alternating, the latter case happening only if $F$ has characteristic two;

    \item if $F$ has characteristic different from two and $a=-2$, then $b$ is similar to $b_{-2}$ if $b$ is non-alternating, and to $h$ if $b$ is alternating;

    \item if $a\neq \pm 2$, then $b$ is similar to $b_a$.
\end{enumerate}
\end{coro}

\begin{proof}
Assume that $a=2$, so that $\chi_{a_b}=(X-1)^2$, and thus $\mu_{a_b}=X-1$ or $(X-1)^2$. If $\mu_{a_b}=(X-1)^2$, then $b$ cannot be symmetric nor alternating, since otherwise, we would have $a_b=\pm \Id_V$ (because alternating forms are skew-symmetric) and $\mu_{a_b}=X\pm 1$. 
But item $(3)$ of the previous lemma then yields a contradiction. Hence $\mu_{a_b}=X-1$, and $b$ is symmetric. The last part of $(1)$ then comes from standard results.

Assume that $F$ has characteristic different from two and $a=-2$, so that $\chi_{a_b}=(X+1)^2$, and thus $\mu_{a_b}=X+1$ or $(X+1)^2$. If $\mu_{a_b}=X+1$, then $b$ is skew-symmetric, and thus isomorphic to $h$. If $\mu_{a_b}=(X+1)^2$, then once again, $b$ cannot be symmetric nor alternating, and we conclude by the previous lemma that $b$ is similar to $b_{-2}$.

Finally, if $a\neq \pm 2$, then $\chi_{a_b}=X^2-aX+1$ is either irreducible or split with simple roots.  In both cases, we have $\mu_{a_b}=\chi_{a_b}$ and $b$ is neither symmetric nor alternating. By the previous lemma, $b$ is similar to $b_a$.
\end{proof}

\begin{lem}\label{chidec}
Let $(V_i,\varphi_i)$ be two non-degenerate bilinear spaces of dimension two, and let $(V,b)=(V_1,\varphi_1)\otimes_F (V_2,\varphi_2)$.

Write $\chi_{a_{\varphi_i}}=X^2-a_iX+1$. Then \[ \chi_{a_b}=X^4-a_1a_2X^3+(a_1^2+a_2^2-2)X^2-a_1a_2+1. \]    
\end{lem}

\begin{proof}
 Let $\alpha_i,\alpha_i^{-1}$ be the (not necessarily distinct) roots  of $\chi_{a_{\varphi_i}}$ in $\overline{F}$. We then have $\alpha_i+\alpha_i^{-1}=a_i$.
 Since $a_b=a_{\varphi_1}\otimes a_{\varphi_2}$, we get \[ \chi_{a_b}=(X-\alpha_1\alpha_2)(X-\alpha_1^{-1}\alpha_2)(X-\alpha_1\alpha_2^{-1})(X-\alpha_1^{-1}\alpha_2^{-1}).  \]
 Thus, \[ \chi_{a_b}=(X^2-a_1\alpha_2X+\alpha_2^2)(X^2-a_1\alpha_2^{-1}X+\alpha_2^{-2}). \]
Since $\alpha_2+\alpha_2^{-1}=a_2$ and $\alpha_2^2+\alpha_2^{-2}=(\alpha_2+\alpha_2^{-1})^2-2=a_2^2-2$, expanding everything yields the desired result.
\end{proof}

We now proceed to show that if $(\mathscr{L}(V'),\sigma_{b'})$ is decomposable, then $(V',b')$ has trivial determinant and decomposable asymmetry.
We need two lemmas.

\begin{lem}\label{quat}
Let $F$ be a field of characteristic different from two, and let $(Q,\tau)$ be a quaternion $F$-algebra with an $F$-linear anti-automorphism.
Then $(Q,\tau)_{\overline{F}}\simeq (\mathscr{L}(V),\sigma_b)_{\overline{F}}$, for some non-degenerate $F$-bilinear space $(V,b)$ of dimension two.
\end{lem}

\begin{proof}
Let $\gamma$ be the canonical involution on $Q$. Then $\tau\circ\gamma^{-1}$ is an automorphism of $Q$, so we may write $\tau={\rm Int}(v)^{-1}\circ \gamma$, for some $v\in Q^\times$. Let $\varphi:Q_{\overline{F}}\overset{\sim}{\to }\mathrm{M}_2(\overline{F})$ be an isomorphism of $\overline{F}$-algebras.

Then $\varphi\circ (\gamma_{\overline{F}})\circ\varphi^{-1}$ is a symplectic involution on $\mathrm{M}_2(\overline{F})$, hence equal to ${\rm Int}(H)^{-1}\circ{}^t$, where $H=\begin{pmatrix}\phantom{-}0 & 1\cr -1 & 0\end{pmatrix}$. It follows that we have  
$$\varphi\circ\tau_{\overline{F}}\circ \varphi^{-1}={\rm Int}(\varphi(v\otimes 1))^{-1}\circ {\rm Int}(H)^{-1}\circ {}^t.$$
Hence, $\varphi$ induces an isomorphism $(Q,\tau)_{\overline{F}}\simeq (\mathscr{L}(\overline{F}^2),\sigma_c)$, where $c:\overline{F}^2\times\overline{F}^2\to \overline{F}$ is the bilinear form whose matrix representation with respect to the canonical basis of $\overline{F}^2$ is $C=H\varphi(v\otimes 1)$.

If $c$ is symmetric or alternating, then $c$ is hyperbolic, hence defined over $F$. Assume now that $c$ is not symmetric nor alternating.
Taking into account that $H^t=-H$, we have $$C^{-1}C^t=\varphi(v\otimes 1)^{-1}H^{-1}\varphi(v\otimes 1)^t H^t=-\varphi(v\otimes 1)^{-1}({\rm Int}(H)^{-1}\circ{}^t)(\varphi(v\otimes 1)),$$
that is $C^{-1}C^t=-\varphi(v\otimes 1)^{-1}\varphi(\gamma_{\overline{F}}(v\otimes 1))=\varphi((-v^{-1}\gamma(v))\otimes 1)$. It follows that the characteristic polynomial of $ C^{-1}C^t$ is the reduced characteristic polynomial of some element of $Q$, hence is an element of $F[X]$. In other words, the characteristic polynomial of $a_c$ has the form $X^2-aX+1$ for some $a\in F$. Since $c$ is not symmetric nor alternating, Corollary \ref{coroba} shows that $c$ is similar to the $F$-bilinear form $b_a$ over $\overline{F}$.
Thus, $(Q,\tau)_{\overline{F}}\simeq (\mathscr{L}(F^2),\sigma_{b_a})_{\overline{F}}$, and we are done. 
\end{proof}

\begin{lem}\label{symfbar}
Let $F$ be a field of characteristic different from two. Then, two non-degenerate $F$-bilinear spaces  are isomorphic/similar over $\overline{F}$ if and only they have conjugate asymmetries over $F$. 
\end{lem}

\begin{proof}
Let $(V,b)$ and $(V',b')$ two non-degenerate $F$-bilinear spaces. Since we have  $a_{b_{\overline{F}}}=(a_b)_{\overline{F}}$ and $a_{b'_{\overline{F}}}=(a_{b'})_{\overline{F}}$, $(V,b)$ and $(V',b')$ have conjugate asymmetries over $F$ if and only if $(V,b)$ and $(V',b')$ have conjugate asymmetries over $\overline{F}$.
  By  \cite[Theorem 8.2 (i)]{Sze}, this last condition is equivalent to say that $(V,b)$ and $(V',b')$ are isomorphic over $\overline{F}$, since $F$ has characteristic different from two. Since similarity and isomorphism coincide over an algebraically closed field, we are done.
\end{proof}

\begin{rem}
This lemma is not true if $F$ has characteristic two. Indeed, symmetric alternating spaces and symmetric non-alternating spaces may never be isomorphic over any field of characteristic two, while they have equal asymmetries.    
\end{rem}

The following corollary was used in the introduction to rephrase our main problem in the split case.

\begin{coro}
Let $F$ be a field of characteristic different from two, and let $(V',b')$ be a non-degenerate $F$-bilinear space. Then, there exists a decomposable $F$-bilinear space $(V,b)$ such that $(\mathscr{L}(V'),\sigma_{b'})_{\overline{F}}\simeq (\mathscr{L}(V),\sigma_{b})_{\overline{F}}$   
if and only if $(V',b')$ has a decomposable asymmetry.
\end{coro}

\begin{proof}
 The condition stated in the lemma is equivalent to ask for $(V',b')$ to be similar to a decomposable $F$-bilinear space over $\overline{F}$. Now, apply Lemma \ref{symfbar} to conclude.
\end{proof}

We are now ready to prove the final result of this section, which was already announced in the introduction.

\begin{prop}\label{adjclosed}
Let $F$ be a field of characteristic different from two, and let $(V',b')$ be a non-degenerate bilinear space of dimension four.
If $(\mathscr{L}(V'),\sigma_{b'})$, resp. $(V',b')$, is decomposable, then $(V',b')$ has a decomposable asymmetry and trivial determinant.
\end{prop}

\begin{proof}
The result is immediate at the level of bilinear spaces.

Assume that $(\mathscr{L}(V'),\sigma_{b'})$ is decomposable, and write $$(\mathscr{L}(V'),\sigma_{b'})\simeq (Q_1,\sigma_1)\otimes_F(Q_2,\sigma_2),$$
where $Q_1,Q_2$ are quaternion $F$-algebras.

By Lemma \ref{quat}, $(Q_i,\sigma_i)_{\overline{F}}\simeq (\mathscr{L}(V_i),\sigma_{\varphi_i})_{\overline{F}}$ for some non-degenerate $F$-bilinear space $(V_i,\varphi_i)$ of dimension two. It follows that $(V',b')$ and $(V_1,\varphi_1)\otimes_F (V_2,\varphi_2)$ are similar over $\overline{F}$. By Lemma \ref{symfbar}, $b'$ and $\varphi_1\otimes \varphi_2$ have conjugate asymmetries. Thus, $(V',b')$ has a decomposable asymmetry.

To prove that $\det(b')$ is trivial, pick an $F$-linear involution $\tau_i$ on $Q_i$, and write $\sigma_i={\rm Int}(v_i)\circ \tau_i$, where $v_i\in Q_i^\times$. We then have $$(\mathscr{L}(V'),\sigma_{b'})\simeq (Q_1\otimes_F Q_2, {\rm Int}(v_1\otimes v_2)\circ (\tau_1\otimes \tau_2)).$$
By \cite[Proposition 8]{CT}, and taking into account that the determinant is invariant under isomorphism, we get 
$$\det(\sigma_{b'})={\rm Nrd}_{Q_1\otimes_FQ_2}(v_1\otimes v_2)\det(\tau_1\otimes \tau_2)\in F^\times/F^{\times 2}.$$
Now, $\tau_1\otimes \tau_2$ is a decomposable $F$-linear involution, so it has trivial determinant. Finally, we get $$\det(\sigma_{b'})={\rm Nrd}_{Q_1\otimes_FQ_2}(v_1\otimes v_2)={\rm Nrd}_{Q_1}(v_1)^2{\rm Nrd}_{Q_2}(v_2)^2=1\in F^\times/F^{\times 2}.$$
Since $\det(b')=\det(\sigma_{b'})$, we are done.
\end{proof}

As explained in the introduction, we would like to know if the converse of the previous result is true: if $(V',b')$ is a non-degenerate $F$-bilinear of dimension four with a decomposable asymmetry and trivial determinant, is $(\mathscr{L}(V'),\sigma_{b'})$, resp. $(V',b')$,  decomposable over $F$ ?

In order to answer this question, we now relate the set of similarity classes of bilinear spaces with a prescribed asymmetry to the set of invertible symmetric elements of some algebra with involution.
  
\section{A useful correspondence}\label{sec-corres}
In this section, $F$ is a field of arbitrary characteristic.

Let $(V,b)$ be a non-degenerate $F$-bilinear space.
As recalled in the introduction, there is a unique endomorphism $a_b$ of $V$, called the \textit{asymmetry} of $b$, such that \[ b(x,y)=b(y, a_b(x)) \ \mbox{ for all }x,y\in V. \]

If $\mathscr{B}$ is an $F$-basis of $V$ and $\mathrm{Mat}(b;\mathscr{B})=B$, then $\mathrm{Mat}(a_b;\mathscr{B})=B^{-1}B^t$.

If $f\in\mathscr{L}(V)$ is an endomorphism of $V$, there is a unique endomorphism $\sigma_b(f)$ such that \[ b(f(x),y)=b(x,\sigma_b(f)(y)) \ \mbox{ for all }x,y\in V . \]

The map $\sigma_b$ is an anti-automorphism of the $F$-algebra $\mathscr{L}(V)$.

It is easy to check that $\sigma_b(a_b)=a_b^{-1}$ and that $\sigma_b^2=\mathrm{Int}(a_b)$. In particular, the restriction of $\sigma_b$ to the centralizer $\mathscr{C}_b=\mathrm{Cent}(a_b)$ of $a_b$ is an involution $*$, which satisfies $a_b^*=a_b^{-1}$.

Thus, for all $f\in \mathscr{C}_b$, we have $b(f(x),y)=b(x, f^*(y))$ for all $x,y\in V$.

Denote by $\mathrm{Sym}(\mathscr{C}_b,*)^\times$ the set of invertible elements of $\mathscr{C}_b$ which are $*$-symmetric.
We now define an equivalence relation on $\mathrm{Sym}(\mathscr{C}_b,*)^\times$.

If $u_1,u_2\in \mathrm{Sym}(\mathscr{C}_b,*)^\times$, we say that $u_1\sim u_2$ 
   if there exist $v\in\mathscr{C}_b^\times$ and $\lambda\in F^\times$ such that $u_2=\lambda vu_1v^*$.

The goal of this section is to relate the corresponding quotient set to the set of similarity classes of certain bilinear spaces. We start with a lemma.

\begin{lem}\label{bu}
Let $(V,b)$ be a non-degenerate bilinear space. For all $u\in \mathrm{Sym}(\mathscr{C}_b,*)^\times$, the map $b_u:V\times V\to F$ defined by \[ b_u(x,y)=b(x,u^{-1}(y)) \ \mbox{ for all }x,y\in V \]
is a non-degenerate bilinear form such that $a_{b_u}=a_b$.
\end{lem}

\begin{proof}
 The bilinearity of $b_u$ and its non-degeneracy are straightforward to establish.
 Moreover, for all $x,y\in V$, we have 
 \[ b_u(y,x)=b(y,u^{-1}(x))=b(u^{-1}(x),a_b(y))=b(x,(u^{-1})^*(a_b(y))).  \]
 Since $u^*=u$, we get \[ b_u(y,x)=b(x,u^{-1}(a_b(y)))=b_u(x,a_b(y)), \]
 and the definition of $a_{b_u}$ yields immediately the equality $a_{b_u}=a_b$.
\end{proof}

\begin{lem}\label{uc}
Let $(V,b)$ be a non-degenerate bilinear space, and let $(V, c)$ a non-degenerate bilinear space such that $a_c=a_b$. Then, there exists a unique $u_c\in \mathrm{GL}(V)$ such that \[ c(x,y)=b(x,u_c^{-1}(y)) \ \mbox{ for all }x,y\in V. \]
Moreover, $u_c\in \mathrm{Sym}(\mathscr{C}_b,*)^\times$.
\end{lem}

\begin{proof}
The uniqueness of $u_c$ comes from the non-degeneracy of $b$. To prove its existence,
 let $\mathscr{B}$ be an $F$-basis of $V$, and let $B$ and $C$ be the matrix representations of $b$ and $c$ with respect to $\mathscr{B}$. Let $u_c\in \mathrm{GL}(V)$ the unique endomorphism of $V$ such that $\mathrm{Mat}(u_c;\mathscr{B})=C^{-1}B$. Then we have \[ c(x,y)=b(x,u_c^{-1}(y)) \ \mbox{ for all }x,y\in V. \]  

 It remains to prove that $u_c$ is an invertible element of $\mathrm{Sym}(\mathscr{C}_b,*)$.
 Note first that for all $f\in\mathscr{L}(V)$, and all $x,y\in V$, we have 
 \[ c(f(x),y)=c(x,\sigma_c(f)(y))=b(x,u_c^{-1}\sigma_c(f)(y)) \] on the one hand, as well as
\[  c(f(x),y)=b(f(x),u_c^{-1}(y))=b(x,\sigma_b(f)u_c^{-1}(y)) \] on the other hand.
Since $b$ is non-degenerate, we get $u_c^{-1}\sigma_c(f)=\sigma_b(f)u_c^{-1}$, that is 
 $\sigma_c=\mathrm{Int}(u_c)\circ \sigma_b$.
Now, we have \[ a_b^{-1}=a_c^{-1}=\sigma_c(a_c)=\sigma_c(a_b)=u_ca_b^{-1}u_c^{-1}.  \]
Consequently, $u_c$ commute with $a_b$. We still have to show that $u_c^*=u_c$.
For all $x,y\in V$, we have \[ b( x, (u_c^*)^{-1}(y))=b(x,(u_c^{-1})^*(y))=b(u_c^{-1}(x),y)=b(y, a_b(u_c^{-1}(x))). \]
Since $u_c$ commutes with $a_b$, and since $a_c=a_b$, for all $x,y\in V$, we get \[ b( x, (u_c^*)^{-1}(y))=b(y, u_c^{-1}(a_b(x)))=c(y,a_b(x))=c(x,y)=b(x,u_c^{-1}(y)). \]
The uniqueness of $u_c$ then implies that $u_c^*=u_c$.
\end{proof}

\begin{thm}\label{corres}
Let $F$ be a field of arbitrary characteristic.

Let $(V,b)$ be a non-degenerate bilinear space. Then, there is a one-to-one correspondence between the set of similarity classes of non-degenerate bilinear spaces $(V',b')$ such that $a_{b'}$ is conjugate to $a_b$, and the quotient set  $\mathrm{Sym}(\mathscr{C}_b,*)^\times/\sim$.

This correspondence sends the equivalence class of $u\in \mathrm{Sym}(\mathscr{C}_b,*)^\times$ to the similarity  class of $(V,b_u)$, and sends the similarity class of $(V',b')$ to the equivalence class of the unique automorphism $u_{b'}$ of $V$ defined by \[  b'(f(x),f(y))=b(x,u_{b'}^{-1}(y)) \ \mbox{ for all }x,y\in V, \]
where $f:V\overset{\sim}{\longrightarrow} V'$ satisfies $a_{b'}=f a_b f^{-1}$.
\end{thm}

\begin{proof}
If $u\in \mathrm{Sym}(\mathscr{C}_b,*)^\times$, we need to show that $b_u$ is a non-degenerate bilinear form such that $a_{b_u}$ and $a_b$ are conjugate, and whose similarity class only depends on the equivalence class of $u$.

Lemma \ref{bu} shows that $b_u$ is non-degenerate, and that $a_{b_u}=a_b$, which is even stronger that the required condition.

Now let $u_1,u_2\in \mathrm{Sym}(\mathscr{C}_b,*)^\times$ such that $u_1\sim u_2$. We are going to show that $b_{u_1}$ and $b_{u_2}$ are similar.

 Let $v\in \mathscr{C}_b^\times$ and $\lambda\in F^\times$ such that $u_2=\lambda vu_1v^*$. Then, for all $x,y\in V$, and taking into account that $(v^*)^{-1}=(v^{-1})^*$, we have \[ \begin{array}{lll}b_{u_2}(v(x),v(y))&=&\lambda^{-1} b(v(x),(v^{-1})^*(u_1^{-1}(y)))\cr &=&\lambda^{-1}b(x,u_1^{-1}(y))\cr &=&\lambda^{-1}b_{u_1}(x,y).\end{array} \]
 Therefore, $b_{u_1}$ and $b_{u_2}$ are similar. 

Now let $(V',b')$ be a non-degenerate bilinear space such that $a_{b'}$ is conjugate to $a_b$. Let $f:V\overset{\sim}{\longrightarrow}V'$ be an isomorphism such that $a_{b'}=f a_bf^{-1}$.

We claim that the bilinear map 
\[ b'_f:(x,y)\in V\times V\mapsto b'(f(x),f(y))\in F \]
is non-degenerate, and that its asymmetry is $a_b$.

The non-degeneracy is easy to establish. Now, for all $x,y\in V$, we have 
\[ \begin{array}{lll}b'_{f}(y,x)&=&b'( f(y),f(x))\cr &=&b'(f(x), a_{b'}(f(y)))\cr &=&b'(f(x),f(a_b(y)))\cr &=&b'_f(x,a_b(y)),\end{array}  \]
and we get $a_{b'_{f}}=a_b$. By Lemma \ref{uc}, there exists an  automorphism $u_{b'}$ of $V$ defined by \[  b'(f(x),f(y))=b(x,u_{b'}^{-1}(y)) \ \mbox{ for all }x,y\in V, \]
and $u_{b'}\in \mathrm{Sym}(\mathscr{C}_b,*)^\times$.

We first prove that the equivalence class of $u_{b'}$ does not depend on the choice of $f$. Let $f_1,f_2:V\overset{\sim}{\longrightarrow}V'$ such that $ a_{b'}=f_1 a_b f_1^{-1}=f_2 a_b f_2^{-1}$. Let us denote by $u_1$ and $u_2$ the elements of $ \mathrm{Sym}(\mathscr{C}_b,*)^\times$ corresponding to the bilinear maps $b'_{f_1}$ and $b'_{f_2}$ respectively.
Set $v=f_2^{-1} f_1$. Then, $v\in\mathscr{C}_b^\times$. Moreover, $f_2=f_1 v^{-1}$, and thus, for all $x,y\in V$, taking into account that $(v^{-1})^*=(v^*)^{-1}$, we have 
\[ \begin{array}{lll}b'_{f_2}(x,y)&=&b'(f_1(v^{-1}(x)),f_1(v^{-1}(y)))\cr &=& b'_{f_1}(v^{-1}(x),v^{-1}(y))\cr  &=&b(v^{-1}(x), u_1^{-1}(v^{-1}(y)))\cr &=&b(x, ((v^*)^{-1}u_1^{-1}v^{-1})(y)), \end{array} \]
and the definition of $u_2$ shows that $u_2=vu_1 v^*$, so in particular $u_1\sim u_2$.

We now prove that the equivalence class of $u_{b'}$ only depends on the similarity class of $(V',b')$. Let $(V'',b'')$ be an $F$-bilinear space which is similar  to $(V',b')$, and let $g: V'\overset{\sim}{\to}V''$ and $\lambda\in F^\times$ such that \[  \lambda b''(g(x'),g(y'))=b'(x',y') \ \mbox{ for all }x',y'\in V'. \]

 Let $f: V\overset{\sim}{\to}V' $ such that $a_{b'}=f a_b f^{-1}$. Then the asymmetry of $b'_{f}$ is $a_b$. Now, we have (with obvious notation) $(\lambda b'')_{gf}=b'_{f}$. It follows from the previous considerations that $u_{\lambda b''}\sim u_{b'}$. But it is clear from the definition that $u_{\lambda b''}=\lambda^{-1} u_{b''}$. Thus, $u_{b''}\sim u_{\lambda b''}\sim u_{b'}$, as required.

It remains to show that the two constructions above are mutually inverse. 

If $u\in \mathrm{Sym}(\mathscr{C}_b,*)^\times$, we need to show that the similarity class of $(V,b_u)$ is sent back to the equivalence class of $u$ via the previous construction. But this follows from the fact $a_{b_u}=a_b$ and the definition of $b_u$.

Now, if $(V',b')$ is a non-degenerate bilinear space such that $a_{b'}$ and $a_b$ are conjugate,  its similarity class is sent to the equivalence class of $u_{b'}$, where $u_{b'}\in  \mathrm{Sym}(\mathscr{C}_b,*)^\times$ satisfies $b'_f=b_{u_{b'}}$, for some well-chosen $f:V\overset{\sim}{\longrightarrow} V'$. Now, the equivalence class of $u_{b'}$ is sent back to the similarity class of $b_{u_{b'}}=b'_f$. But by very definition, $b'_f$ is isomorphic, hence similar, to $b'$.

This concludes the proof.
\end{proof}

We now use the previous theorem to start investating Question 2 in some special cases.
We will assume until the end of this section that $F$ has characteristic different from two.

{\bf Notation. }If $\alpha\in F$, we will denote by $F[\sqrt{\alpha}]$ the quadratic $F$-algebra $F[T]/(T^2-\alpha)$. We will denote by $\sqrt{\alpha}$ the class of $T$, and by $-$ the unique non-trivial automorphism of this quadratic $F$-algebra.

If $(V,\varphi)$ is a non-degenerate bilinear space of dimension two, then Lemma \ref{lemba} and Corollary \ref{coroba} show that $a_\varphi$ is either equal to $\pm \mathrm{Id}_V$, or is a cyclic endomorphism with minimal polynomial equal to $X^2-aX+1$, for some $a\in F\setminus\{2\}$.

If $a_\varphi=\pm \mathrm{Id}_V$, then  $(\mathscr{C}_\varphi,*)=(\mathscr{L}(V),\sigma_\varphi)$. 

If $a_\varphi$ is cyclic, with minimal polynomial equal to $X^2-aX+1$, then $\mathscr{C}_\varphi=F[a_\varphi]$. Moreover, since $a_\varphi^*=a_\varphi^{-1}=a\mathrm{Id}_V-a_\varphi$, we have $(2a_\varphi-a\mathrm{Id}_V)^*=-(2a_\varphi-a\mathrm{Id}_V)$. Since $(2a_\varphi-a\mathrm{Id}_V)^2=(a^2-4)\mathrm{Id}_V$, we then get an isomorphism of $F$-algebras with involutions $$(F[\sqrt{a^2-4}], - )\simeq (\mathscr{C}_\varphi,*)$$
sending $\sqrt{a^2-4}$ to $2a_\varphi-a\mathrm{Id}_V$.
In particular, if $a^2-4$ is a non-zero square, we have $(\mathscr{C}_\varphi, *)\simeq (F^2, s)$, where $s$ is the switch.

Finally, it is clear that, if $(V_1,\varphi_1), (V_2,\varphi_2)$ are two non-degenerate  $F$-bilinear spaces, then $(\mathscr{C}_{\varphi_1\otimes\varphi_2},*)=(\mathscr{C}_{\varphi_1},*)\otimes_F(\mathscr{C}_{\varphi_2},*)$.

We continue with a convenient definition.

\begin{defn}\label{genasym}
We say that a non-degenerate $F$-bilinear space $(V',b')$ of dimension four has {\it a generic decomposable asymmetry} if its asymmetry is conjugate to the asymmetry of $(F^2,b_{a_1})\otimes_F(F^2,b_{a_2})$, for some $a_1,a_2\in F\setminus\{-2,2\}$  such that $\alpha_i\overset{\mathrm{def}}{=}a_i^2-4$ is not a square in $F^\times$ for $i=1,2$.
\end{defn}

We then have the following proposition.

\begin{prop}\label{non-gen}
Let $(V',b')$ be a non-degenerate $F$-bilinear space of dimension four with a non-generic decomposable asymmetry. Then, $(V'b')$ is decomposable if and only if $\det(b')$ is trivial.
\end{prop}

\begin{proof}

Assume that the asymmetry of $(V',b')$ is conjugate to the asymmetry of $(V_1,\varphi_1)\otimes (V_2,\varphi_2)$.
 
If $(V',b')$ is symmetric or skew-symmetric, then $(V',b')$ has non-generic decomposable asymmetry, and the result is known to be true.

Hence, we may assume that one the $\varphi_i$'s is not symmetric nor skew-symmetric.
We may also assume without loss of generality that $V_1=V_2=F^2$. 
Finally, note that any symmetric bilinear space of dimension two has same asymmetry as $(F^2, \langle 1,1\rangle)$, and any skew-symmetric bilinear space
of dimension two has same asymmetry as $(F^2,h)$.
Hence, using Lemma \ref{lemba}, the symmetry of the tensor product of bilinear spaces and the definition of a generic asymmetry, we see that the remaining cases are the following ones:

(i) $(V_1,\varphi_1)\otimes_F(V_2,\varphi_2)=(F^2,\varphi)\otimes (F^2, b_a)$, where $a\in F\setminus\{-2, 2\}$ is such that $\alpha=a^2-4$ is a square, and $(F^2,\varphi)$ is an arbitrary $F$-bilinear space of dimension two;

(ii) $(V_1,\varphi_1)\otimes_F(V_2,\varphi_2)=(F^2,\varphi)\otimes (F^2, b_a)$, where $a\in F\setminus\{-2, 2\}$ is such that $\alpha=a^2-4$ is not a square and $\varphi=\langle 1,1\rangle$ or $h$;

(iii) $(V_1,\varphi_1)\otimes_F(V_2,\varphi_2)=(F^2,\varphi)\otimes (F^2, b_{-2})$,  where  $\varphi=\langle 1,1\rangle, h, b_{-2}$ or $b_a$, where 
$a\in F\setminus\{-2,2\}$ is such that $\alpha$ is not a square.

In the sequel, we will denote by $\mathscr{E}=(\varepsilon_1,\varepsilon_2)$  the canonical basis of $F^2$, and by 
$\mathscr{B}_0$ the $F$-basis of $F^2\otimes_FF^2$ defined by \[ \mathscr{B}_0= (\varepsilon_1\otimes \varepsilon_1,\varepsilon_1\otimes \varepsilon_2,\varepsilon_2\otimes \varepsilon_1,\varepsilon_2\otimes \varepsilon_2). \]

Recall that, if $a\in F\setminus\{-2,2\}$, we have  \[ \mathrm{Mat}(a_{b_a}; \mathscr{E})=A_a\overset{\mathrm{def}}{=}\begin{pmatrix}
a-1 & -1 \cr 2-a & \phantom{-}1
\end{pmatrix}. \]

- Case (i). In this case, we have $(\mathscr{C}_{b_a},*)\simeq (F^2, s)$, where $s$ is the switch, and therefore 
\[ (\mathscr{C}_{\varphi_1\otimes\varphi_2},*)\simeq (\mathscr{C}_\varphi^2 ,\rho), \]
where $\rho((x,y))=(y^*,x^*)$ for all $x,y\in \mathscr{C}_\varphi$.

Now, $\mathrm{Sym}(\mathscr{C}_\varphi^2 ,\rho)^\times=\{(x,x^*)\mid x\in\mathscr{C}_\varphi^\times\}$.
But for all $x\in \mathscr{C}_\varphi^\times$, we have $(x,x^*)=(x,1)\rho((x,1))$. Thus, $\mathrm{Sym}(\mathscr{C}_{\varphi_1\otimes \varphi_2},*)^\times/\sim$ is trivial, and we conclude using Theorem \ref{corres} that $b'$ is similar to $\varphi_1\otimes \varphi_2$, hence is decomposable.

- Case (ii). The matrix representation of the asymmetry of $\langle 1,1\rangle\otimes b_a$ with respect to $\mathscr{B}_0$ is $I_2\otimes A_a$, while the
 matrix representation of the asymmetry of $h\otimes b_{-a}$ is $-I_2\otimes A_{-a}$.
 It is then easy to see that $-I_2\otimes A_{-a}=P^{-1}(I_2\otimes A_a)P$, where $P=I_2\otimes\begin{pmatrix}
 1 & \frac{2}{a+2} \cr 0 & \frac{a-2}{a+2}
 \end{pmatrix}$. Hence, replacing $a$ by $-a$, one may assume without loss of generality that $\varphi=\langle 1,1\rangle$.
In this case, we have \[(\mathscr{C}_{\varphi_1\otimes \varphi_2},*)=(\mathscr{L}(F^2),\sigma_{\langle 1,1\rangle})\otimes_F (F[a_{b_a}],*)\simeq (\mathrm{M}_2(F[\sqrt{\alpha}]), -^t),\] where $-$ is the non-trivial $F$-automorphism of $F[\sqrt{\alpha}]$. This isomorphism sends  $f\otimes \mathrm{Id}$ to $\mathrm{Mat}(f;\mathscr{E})$, and $\mathrm{Id}\otimes a_{b_a}$ to  $\dfrac{a+\sqrt{\alpha}}{2}I_2$.

Now, invertible symmetric elements of $(\mathscr{C}_{\varphi_1\otimes \varphi_2},*)$ identify to invertible $2\times 2$ hermitian matrices. 
Note that the symmetric elements of $(F[\sqrt{\alpha}],-)$ are exactly the elements of $F$.

Let $M=\begin{pmatrix}
a & b \cr \overline{b}&  c
\end{pmatrix}\in \mathrm{M}_2(F[\sqrt{\alpha}])$ be an invertible hermitian matrix. In particular, $a,c\in F$.

If $a\in F^\times$, we have $M=P\begin{pmatrix}
a & 0 \cr  0 & c- \overline{b}b a^{-1}
\end{pmatrix} \overline{P}^t$, where $P=\begin{pmatrix}
1 & 0 \cr \overline{b} a^{-1} & 1
\end{pmatrix}$, 
while if $c\in F^\times$, we have $M=P\begin{pmatrix}
a- \overline{b}bc^{-1}& 0 \cr 0 & c
\end{pmatrix}\overline{P}^t$, where $P=\begin{pmatrix}
1 &  bc^{-1} \cr 0 & 1
\end{pmatrix}$.
Finally, if $a=c=0$, then $b\in F[\sqrt{\alpha}]^\times$, and we have $M=P\begin{pmatrix}
\frac{1}{2} & \phantom{-}0 \cr0 & -\frac{1}{2} 
\end{pmatrix}\overline{P}^t$, where $P=\begin{pmatrix}
b & \phantom{-}b \cr 1 & -1
\end{pmatrix}$.

Using the isomorphism above, we see that every symmetric invertible element of $(\mathscr{C}_{\varphi_1\otimes \varphi_2},*)$ is equivalent to $f\otimes \mathrm{Id}$, where $f$ is an automorphism of $F^2$ which is self-adjoint with respect to the unit form. It follows easily from the explicit correspondence that the corresponding similarity class is represented by  $\psi\otimes b_a$, where $(F^2,\psi)$ is a symmetric bilinear space. Therefore, $b'$ is decomposable.

- Case (iii).
Note that $(\mathscr{C}_{b_{-2}},*)$ is isomorphic to $(F[\sqrt{0}],-)$, so that $(\mathscr{C}_{\varphi_1\otimes \varphi_2},*)\simeq (\mathscr{C}_\varphi,*)\otimes_F(F[\sqrt{0}],-)$.

Assume first that $\varphi=\langle 1,1\rangle$. Then the proof given in Case (ii) is still valid for $\alpha=0$, and we are done. Assume now that we are in the three remaining cases. It is easy to check that, in each of these cases, the set of symmetric elements of $(\mathscr{\varphi},*)$ identifies to $F$.
Setting $\varepsilon=\overline{T}$, we get that the symmetric elements of $(\mathscr{C}_\varphi, *)\otimes_F (F[\sqrt{0}],-)$ have the form $a\mathrm{Id}\otimes 1 +b\otimes \varepsilon$, where $a\in F$ and $b^*=-b$. Such an element being invertible if and only if $a\in F^\times$, we finally get that any symmetric invertible element has the form $u=a(\mathrm{Id}\otimes 1+z\otimes \varepsilon)$, where $z^*=-z$ and $a\in F^\times$.
But we have $$a(\mathrm{Id}\otimes 1+\frac{z}{2}\otimes \varepsilon)(\mathrm{Id}\otimes 1+\frac{z}{2}\otimes \varepsilon)^*=a(\mathrm{Id}\otimes 1+\frac{z}{2}\otimes \varepsilon)^2=a(\mathrm{Id}\otimes1 +z\otimes \varepsilon)=u,$$
so the class of $u$ is trivial. We conclude as before that $b'$ is similar to $\varphi\otimes b_{-2}$, hence is decomposable.
\end{proof}

\section{Determinant and decomposability}\label{sec-det-gen}

In all this section, we will assume that $F$ has characteristic different from two.
The main goal here is to answer Questions 1-3. We start with Question 2, and will concentrate on the case of bilinear spaces with a generic decomposable asymmetry, in view of Proposition \ref{non-gen}.

Let $a_1,a_2\in F\setminus\{2,-2\}$, and set $\alpha_i=a_i^2-4$. We assume in this section that $\alpha_1,\alpha_2$ are not squares in $F^\times$.

We will set $L=F[\sqrt{\alpha_1}]\otimes_F F[\sqrt{\alpha_2}]$, as well as $K=F[\sqrt{\alpha_1}\otimes \sqrt{\alpha_2}]$. 

Keeping the notation of the previous sections, set $(V,b)=(F^2,b_{a_1})\otimes (F^2, b_{a_2})$.

If $a\in F\setminus\{-2,2\}$, recall that we have \[ \mathrm{Mat}(a_{b_a}; \mathscr{E})=A_a\overset{\mathrm{def}}{=}\begin{pmatrix}
    a-1 & -1 \cr 2-a & \phantom{-}1
\end{pmatrix},\] as well as $\mathrm{Mat}(a_b;\mathscr{B}_0)=A_{a_1}\otimes A_{a_2}.$

The computations made in the previous section show that
 we have an isomorphism of $F$-algebras with involutions $ (\mathscr{C}_b,*)\simeq (L,\sigma)$,
 where $\sigma$ is the unique $F$-algebra automorphism of $L$ sending $\sqrt{\alpha_i}$ to $-\sqrt{\alpha_i}$ for $i=1,2$.

We then immediately get the following proposition.

\begin{prop}\label{corres4}
Let $L=F[\sqrt{\alpha_1}]\otimes_F F[\sqrt{\alpha_2}]$, and let $K=F[\sqrt{\alpha_1}\otimes \sqrt{\alpha_2}]$.

There is a one-to-one correspondence between the set of similarity classes of non-degenerate bilinear spaces $(V',b')$ such that $a_{b'}$ is conjugate to $a_b$, and the quotient group  $K^\times/F^\times N_{L/K}(L^\times)$.

This correspondence sends the equivalence class of $r+s\sqrt{\alpha_1}\otimes\sqrt{\alpha_2}\in K^\times$ to the similarity  class of $b_u$, where \[ u=r\Id+s (-a_1\Id+2a_{b_{a_1}})\otimes (-a_2\Id+2a_{b_{a_2}})\in\mathrm{Sym}(\mathscr{C}_b,*)^\times, \] and sends the similarity class of $(V',b')$ to the equivalence class of $r+s \sqrt{\alpha_1}\otimes \sqrt{\alpha_2}$, where $r,s\in F$ are such  that  $u_{b'}=r\Id+s (-a_1\Id+2a_{b_{a_1}})\otimes (-a_2\Id+2a_{b_{a_2}})$ is the unique automorphism $u_{b'}$ of $F^2\otimes_FF^2$ defined by \[  b'(f(x),f(y))=b(x,u_{b'}^{-1}(y)) \ \mbox{ for all }x,y\in V, \]
where $f:V\overset{\sim}{\longrightarrow} V'$ satisfies $a_{b'}=f a_b f^{-1}$.
\end{prop}

\begin{rem}\label{a1a2square}
Assume that $\alpha_2\alpha_1^{-1}$ is a non-zero square of $F$, so that we have  $L=F[\sqrt{\alpha_1}]\otimes_FF[\sqrt{\alpha_1}]$ and $K=F[\sqrt{\alpha_1}\otimes \sqrt{\alpha_1}]$. Set $K_1=F[\sqrt{\alpha_1}]$.

Then, there exists a unique isomorphism of $F$-algebras $$\tau:L\overset{\sim}{\to }K_1\times K_1$$ sending $z_1\otimes z_2$ to $(z_1z_2,z_1\overline{z}_2)$, where $-$ is the unique non-trivial $F$-automorphism of $K_1$.

This isomorphism then induces an $F$-algebra isomorphism $K\overset{\sim}{\to} F\times F$, sending $r+s\sqrt{\alpha_1}\otimes\sqrt{\alpha_1}$ to $(r+s\alpha_1,r-s\alpha_1)$. Note that for all $\lambda_1,\lambda_2\in F^\times$, we have $$\tau^{-1}((\lambda_1,\lambda_2))=\frac{\lambda_1+\lambda_2}{2}+ \frac{\lambda_1-\lambda_2}{2\alpha_1}\sqrt{\alpha_1}\otimes\sqrt{\alpha_1}.$$

Moreover, for all $z_1,z_2\in K_1$, we have $\sigma(\tau(z_1\otimes z_2))=(\overline{z_1z_2}, \overline{z_1\overline{z}_2})$.
It follows easily that we have an isomorphism of algebras with involution $$(L,\sigma)\simeq (K_1\times K_1,\rho),$$ where 
$\rho((z_1,z_2))=(\overline{z}_1,\overline{z}_2)$ for all $z_1,z_2\in K_1$.

In particular, we have $N_{L/K}(L^\times)=N_{K_1/F}(K_1^\times)\times N_{K_1/F}(K_1^\times)$.
Thus, $\tau$ induces a group isomorphism $$K^\times /F^\times N_{L/K}(L^\times)\simeq (F^\times \times F^\times)/F^\times(N_{K_1/F}(K_1^\times)\times N_{K_1/F}(K_1^\times)).$$
Now, the group morphism $(\lambda_1,\lambda_2)\in F^\times \times F^\times \mapsto \lambda_1\lambda_2^{-1}\in F^\times$ induces a group isomorphism 
between the latter group and $ F^\times/ N_{K_1/F}(K_1^\times)$,
the inverse map being induced by $\lambda\in F^\times \mapsto (\lambda,1)\in F^\times \times F^\times$.

All in all, we get a group isomorphism $$\Theta:K^\times/F^\times N_{L/K}(L^\times)\overset{\sim}{\to} F^\times/ N_{K_1/F}(K_1^\times),$$ sending the class of $r+s \sqrt{\alpha_1}\otimes \sqrt{\alpha_1}\in K^\times$ to the class of $\dfrac{r+s\alpha_1}{r-s\alpha_1}$.

One easily check that $\Theta^{-1}$ sends the class of $\lambda\in F^\times$ to the class of $\dfrac{\lambda+1}{2}+\dfrac{\lambda-1}{2\alpha_1}\sqrt{\alpha_1}\otimes \sqrt{\alpha_1}$.

Note for later use that for all $\lambda\in F^\times$, we have $N_{K/F}(\Theta^{-1}(\lambda))=\lambda$.
 \end{rem}

We now determine all the decomposable bilinear spaces $(V',b')$ such that $(V,b)$ and $(V',b')$ have conjugate asymmetries.

\begin{lem}\label{decasym}
The bilinear forms $b$ and $b_{-a_1}\otimes b_{-a_2}$ have conjugate asymmetries.

Moreover, if $(V',b')$ is a non-degenerate bilinear space such that $(V,b)$ and $(V',b')$ have conjugate asymmetries, then $b'$ is decomposable if and only if it is similar to $b$ or to $b_{-a_1}\otimes b_{-a_2}$.
\end{lem}

\begin{proof}
We have $\mathrm{Mat}(b_{-a_1}\otimes b_{-a_2};\mathscr{B}_0)=A_{-a_1}\otimes A_{-a_2}$.

Let \[ P=\begin{pmatrix}
1&  -\frac{2}{-2+a_2} & -\frac{2}{-2+a_1} &  \frac{4}{(-2+a_1)(-2+a_2)}\cr 0&  \frac{2+a_2}{-2+a_2} &  0&  -\frac{4+2a_2}{(-2+a_1)(-2+a_2)}\cr 0& 0&  \frac{2+a_1}{-2+a_1} &  -\frac{4+2a_1}{(-2+a_1)(-2+a_2)}\cr  0 &  0 &  0 & \frac{(2+a_1)(2+a_2)}{(-2+a_1)(-2+a_2)}\end{pmatrix}. \]

Since $a_1,a_2\neq \pm 2$, the matrix $P$ is well-defined and invertible. Now, one check  that \[  A_{-a_1}\otimes A_{-a_2}=P(A_{a_1}\otimes A_{a_2})P^{-1}, \]
which proves the first part of the lemma.

Since the characteristic polynomial of $a_{b_{a_i}}$ is $X^2-a_iX+1$, the characteristic polynomial of $a_b$ is \[ \chi_{a_b}=X^4-a_1a_2X^3+(a_1^2+a_2^2-2)X^2-a_1a_2+1, \]
by Lemma \ref{chidec}. By assumption, we have $b'\simeq \varphi'_1\otimes \varphi'_2$, where $(V'_i,\varphi'_i)$ is a non-degenerate bilinear space of dimension two. If $\chi_{a_{\varphi'_i}}=X^2-a'_iX+1$, then, by loc.cit., we have \[ \chi_{a_{b'}}=X^4-a'_1a'_2X^3+({a'_1}^2+{a'_2}^2-2)X^2-a'_1a'_2+1. \]
Since $a_b$ and $a_{b'}$ are supposed to be conjugate, they have same characteristic polynomial, and we get $a'_1a'_2=a_1a_2$ and ${a'_1}^2+{a'_2}^2=a_1^2+a_2^2$.

 We then have $(a'_1+a'_2)^2=(a_1+a_2)^2$
and $(a'_1-a'_2)^2=(a_1-a_2)^2$, and thus $a'_1+a'_2=\varepsilon(a_1+a_2)$ and $a'_1-a'_2=\varepsilon'(a_1-a_2)$, where $\varepsilon,\varepsilon'=\pm 1$.

If $\varepsilon'=\varepsilon$, we then get $(a'_1,a'_2)=(\varepsilon a_1,\varepsilon a_2) $, and if $\varepsilon'=-\varepsilon$, we get $(a'_1,a'_2)=(\varepsilon a_2,\varepsilon a_1)$.

 Since $a_1,a_2\neq \pm 2$, Corollary \ref{coroba} implies the desired result, taking into account that tensor product of bilinear spaces is commutative up to isomorphism or similarity.
 \end{proof}

Using Proposition \ref{corres4}, we derive the following result.

\begin{prop}\label{propdecomp}
 Let $(V',b')$ be a non-degenerate bilinear space such that $(V,b)$ and $(V',b')$ have conjugate asymmetries. Then $b'$ is decomposable if and only if its similarity class corresponds to the trivial class or the class of $\sqrt{\alpha_1}\otimes\sqrt{\alpha_2}$ in the quotient group $K^\times/F^\times N_{L/K}(L^\times)$.
\end{prop}

\begin{proof}
 Using Lemma \ref{decasym}, we see that it remains to check that the similarity class of $b_{-a_1}\otimes b_{-a_2}$ corresponds to the class of  $\sqrt{\alpha_1}\otimes\sqrt{\alpha_2}$.

 For, we follow the description of the correspondence given in Proposition \ref{corres4}.

 Let $B=\mathrm{Mat}(b;\mathscr{B}_0)$ and let $B'=\mathrm{Mat}(b_{-a_1}\otimes b_{-a_2};\mathscr{B}_0)$.
 Let $P$ be the invertible matrix given in the proof of Lemma \ref{decasym}, and let  $u\in \mathrm{Sym}(\mathscr{C}_b,*)^\times$ be the invertible symmetric element whose class corresponds to the similarity class of $b_{-a_1}\otimes b_{-a_2}$. Then $U=\mathrm{Mat}(u;\mathscr{B}_0)$ satisfies  $P^t B' P=BU^{-1}$. 
 
 Direct computations show that \[ U=\dfrac{1}{(2+a_1)(2+a_2)}(-a_1I_2+2A_{a_1})\otimes (-a_2I_2+2A_2). \]
 In other words, $u=\dfrac{1}{(2+a_1)(2+a_2)}(-a_1\Id+2a_{b_{a_1}})\otimes (-a_2\Id+2a_{b_{a_2}})$.
 By Proposition \ref{corres4}, the similarity class of $b'$ then corresponds to the class of \[\dfrac{1}{(2+a_1)(2+a_2)}\sqrt{\alpha_1}\otimes\sqrt{\alpha_2}\] in $K^\times/F^\times N_{L/K}(L^\times)$, which is also the class of  $\sqrt{\alpha_1}\otimes\sqrt{\alpha_2}$.
This yields the desired result.
\end{proof}

Our next goal is to show that having a decomposable asymmetry and a trivial determinant is not enough to be decomposable.
We will need a tractable way to determine whether a class in $K^\times/F^\times N_{L/K}(L^\times)$ is trivial or not.
The following lemma will be helpful.

\begin{lem}\label{biquad}
Let $u\in K^\times$. Then $u\in F^\times N_{L/K}(L^\times)$ if and only if $N_{K/F}(u)\in N_{K_1/F}(K_1^\times)$, where $K_1=F[\sqrt{\alpha_1}]$.
\end{lem}

\begin{proof}
Let $M=F[\sqrt{\alpha_1}\otimes 1]$. Any element of $M$ has the form $y\otimes 1$, for a unique $y\in K_1$.
Since we have \[ N_{K_1/F}(y)=N_{F[\sqrt{\alpha_1}\otimes 1]/F}(y\otimes 1),\]
 the statement we would like to prove is equivalent to the following one :  for all $u\in K^\times$, we have $u\in F^\times N_{L/K}(L^\times)$ if and only if $N_{K/F}(u)\in N_{M/F}(M^\times)$. 

Assume first that $u=\lambda N_{L/K}(x)$, where $\lambda\in F^\times$ and $x\in L^\times$.
Then, we have \[ N_{K/F}(u)=\lambda^2 N_{L/F}(x)=\lambda^2N_{M/F}(N_{L/M}(x))=N_{M/F}(\lambda N_{L/M}(x)). \]

We now prove the reverse implication. Let $u\in K^\times$ such that $N_{K/F}(u)=N_{M/F}(v)$ for some $v\in M^\times$.

 Denote by $\sigma_1,\sigma_2,\sigma_3$ the $F$-automorphisms of $L$ which are uniquely defined by the equalities \[ \sigma_1(\sqrt{\alpha_1}\otimes 1)=\sqrt{\alpha_1}\otimes 1, \ \sigma_1(1\otimes \sqrt{\alpha_2})=-1\otimes \sqrt{\alpha_2}, \]   
 \[ \sigma_2(\sqrt{\alpha_1}\otimes 1)=-\sqrt{\alpha_1}\otimes 1, \ \sigma_2(1\otimes \sqrt{\alpha_2})=1\otimes \sqrt{\alpha_2}, \]   
 \[ \sigma_3(\sqrt{\alpha_1}\otimes 1)=-\sqrt{\alpha_1}\otimes 1, \ \sigma_3(1\otimes \sqrt{\alpha_2})=-1\otimes\sqrt{\alpha_2}. \]   

We then get 
\[  N_{K/F}(u)=u\sigma_2(u)=N_{M/F}(v)=v\sigma_2(v),  \]
that is $uv^{-1}\sigma_2(uv^{-1})=1$.
Now, $\sigma_2$ is the unique non-trivial 
automorphism of the quadratic \'{e}tale $F[1\otimes \sqrt{\alpha_2}]$-algebra $L$, so the previous equality may be rewritten as $N_{L/F[1\otimes\sqrt{\alpha_2}]}(uv^{-1})=1$. By Hilbert 90, there exists $w\in L^\times$ such that $uv^{-1}=w \sigma_2(w)^{-1}$.
We then have $u=vw \sigma_2(w)^{-1}$. Since $u\in K^\times$, we have $u=\sigma_3(u)$, and therefore \[  vw\sigma_2(w)^{-1}=\sigma_3(v)\sigma_3(w)\sigma_1(w)^{-1}, \]
that is, $vw \sigma_1(w)=\sigma_3(v)\sigma_3(w)\sigma_2(w)=\sigma_3(vw\sigma_1(w))$. 

Note that $w\sigma_1(w)$ is fixed by $\sigma_1$, so $w\sigma_1(w)\in M^\times$, and thus $vw\sigma_1(w)$ is an element of $M^\times$ which is fixed by $\sigma_3$.
It readily follows that $vw\sigma_1(w)=\lambda\in F^\times$.
Hence, we have $u=v w \sigma_2(w)^{-1}=\lambda(\sigma_1(w)\sigma_2(w))^{-1}$.
Setting $x=\sigma_1(w)^{-1}$, we then obtain $u=\lambda x \sigma_3(x)=\lambda N_{L/K}(x)$ as required, since $\sigma_3$ is the unique non-trivial automorphism of $L/K$.
\end{proof}

The standard properties of quaternion algebras then yield the following result.

\begin{prop}\label{critriv}
 Let $u\in K^\times$. Then, the class of $u$ in $K^\times/F^\times N_{L/K}(L^\times)$ is trivial if and only if the quaternion $F$-algebra $(\alpha_1,N_{K/F}(u))$ splits.
\end{prop}

\begin{rem}\label{remdec}
The previous results and properties of quaternion algebras  show that $b_{a_1}\otimes b_{a_2}$ and $b_{-a_1}\otimes b_{-a_2}$  are similar if and only if the quaternion $F$-algebra $(\alpha_1,\alpha_2)$ splits.
 \end{rem}

Before stating the next result of this section, we need another lemma.

\begin{lem}\label{lemdet}
     Let $(V',b')$ be a non-degenerate bilinear space such that $(V,b)$ and $(V',b')$ have conjugate asymmetries. Then $b'$ has trivial determinant.
\end{lem}

\begin{proof}
 If $x+y\sqrt{\alpha_1}\otimes \sqrt{\alpha_2}\in K^\times$, it corresponds to an endomorphism $u$ of $F^2\otimes_F F^2$ such that \[ \mathrm{Mat}(u;\mathscr{B}_0)=xI_4+y(-a_1I_2+2A_{a_1})\otimes (-a_2I_2+2A_{a_2})\in\mathscr{C}_b. \] A direct computation shows that this matrix has square determinant.
Now, the matrix of $b_u$ with respect to this same basis is simply $BU^{-1}$ (with obvious notation), which has square determinant, since $B$ and $U$ have square determinant. This concludes the proof.
\end{proof}

We then get the following theorem.

\begin{thm}\label{thmnondec}
 If $F$ is a number field, there are infinitely many similarity classes of non-degenerate bilinear spaces $(V',b')$ such that $a_b$ and $a_{b'}$ are conjugate and $\det(b')$ is trivial.

Moreover, the only decomposable ones are $b$ and $b_{-a_1}\otimes b_{-a_2}$, up to similarity.
\end{thm}

\begin{proof}
By Proposition \ref{corres4}, we need to prove that  $K^\times/F^\times N_{L/K}(L^\times)$ is infinite.   

We start with the most difficult case, which the case where $L$ is a field.

 \textbf{Claim. }There exists an infinite set $T$ of prime ideals of $\mathscr{O}_F$ such that :
 
 \begin{enumerate}
\item for all $\mathfrak{p}\in T$, $\mathfrak{p}$ stays inert in $F[\sqrt{\alpha_1}]$

\item  for all $\mathfrak{p}\in T$, $\mathfrak{pp}$ is non-dyadic and $v_\mathfrak{p}(\alpha_1)=0$ 

\item  for all $\mathfrak{p}\in T, \mathfrak{p}$ splits in $K$

\item  there exists a family $(u_\mathfrak{p})_{\mathfrak{p}\in T}$ of elements of $K^\times$ such that, for all $\mathfrak{p},\mathfrak{q}\in T$, we have \[ v_\mathfrak{q}(N_{K/F}(u_\mathfrak{p}))=\Bigl\lbrace\begin{array}{ccc}1 & \mbox{ if }&\mathfrak{q}=\mathfrak{p}\cr 0 & \mbox{ if } & \mathfrak{q}\neq\mathfrak{p}. \end{array} \]
 \end{enumerate}

Assume the claim is proved. 

We now proceed to show that the various classes of the elements $u_\mathfrak{p}, \ \mathfrak{p}\in T$ are non-trivial and pairwise distinct, which will be enough to conclude, since $T$ is infinite.

 Let $\mathfrak{p}\in T$. Conditions 1. and 4. imply that  the residue of $(\alpha_1, N_{K/F}(u_\mathfrak{p}))$ at $\mathfrak{p}$ is the square class of $\alpha_1$ in $(\mathscr{O}_{F_\mathfrak{p}}/\mathfrak{p}_\mathfrak{p})^\times$. This class cannot be trivial, since otherwise the residual degree of $F[\sqrt{\alpha_1}]/F$ at $\mathfrak{p}$ would be $1$, which would contradict condition 1.
In particular, $(\alpha_1, N_{K/F}(u_\mathfrak{p}))$ does not split. 
 Hence, the class of $u_\mathfrak{p}$ in $K^\times/F^\times N_{L/K}(L^\times)$ is non-trivial by Proposition \ref{critriv}.

We now proceed to show that for any pair of distinct prime ideals $\mathfrak{p},\mathfrak{q}\in T$, $u_\mathfrak{p}$ and $u_\mathfrak{q}$ have distinct classes in the quotient group.
We thus need to show that $(\alpha_1, N_{K/F}(u_\mathfrak{p} u^{-1}_\mathfrak{q}))$ is not split. 
But we have \[  v_\mathfrak{p}(N_{K/F}(u_\mathfrak{p} u^{-1}_\mathfrak{q}))= 1, \] since $\mathfrak{p}\neq\mathfrak{q}$.
Using similar arguments as before, we reach the desired conclusion.
The second part has been proven in Lemma \ref{decasym}, and we are done.

It remains to prove the claim. 

Since $\mathrm{Gal}(L/F)$ is canonically isomorphic to $\mathrm{Gal}(F[\sqrt{\alpha_1}]/F)\times \mathrm{Gal}(K/F)$, Chebotarev's density theorem shows that the set $S$ of non-zero prime ideals $\mathfrak{p}$ of $\mathscr{O}_F$ satisfying $1-3$ is infinite.

 If $\mathfrak{p}\in S$, let $\mathfrak{P}_\mathfrak{p}$ one of the two distinct prime ideals of $\mathscr{O}_K$ lying above $\mathfrak{p}$. Since the class group of $K$ is finite, there exist $\mathfrak{p}_1,\ldots\mathfrak{p}_r\in S$  such that, for all $\mathfrak{p}\in S$, there exist $u_\mathfrak{p}\in K^\times$ and $i_\mathfrak{p}\in\bigl\llbracket 1,r\bigr\rrbracket$ satisfying $\mathfrak{P}_\mathfrak{p}=u_\mathfrak{p} \mathfrak{P}_{\mathfrak{p}_{i_\mathfrak{p}}}$.
Taking relative norms, we get \[ \mathfrak{p}=N_{K/F}(u_\mathfrak{p})\mathfrak{p}_{i_\mathfrak{p}}. \]

 Let $T=S\setminus\{\mathfrak{p}_1,\ldots,\mathfrak{p}_r\}$.
 The previous equality then shows that, for all $\mathfrak{p},\mathfrak{q}\in T$, we have \[ v_\mathfrak{q}(N_{K/F}(u_\mathfrak{p}))=\Bigl\lbrace\begin{array}{ccc}1 & \mbox{ if }&\mathfrak{q}=\mathfrak{p}\cr 0 & \mbox{ if } & \mathfrak{q}\neq\mathfrak{p} \end{array}, \]
  and we are done.
  
  It remains to deal with the case where $L$ is not a field. Since $\alpha_1,\alpha_2\in F^\times$ are not squares by assumptions, this is equivalent to say that $\alpha_2\alpha_1^{-1}$ is a non-zero square. In this case, Remark \ref{a1a2square} implies that the desired result boils down to prove that the group $F^\times /N_{K_1/F}(K_1^\times)$ is infinite, where $K_1=F[\sqrt{\alpha_1}]$.
  
Note that the class of $\lambda\in F^\times$ is trivial in the latter quotient group if and only if the quaternion $F$-algebra $(\alpha_1,\lambda)$ splits. 

Let $T$ be the set of prime ideals $\mathfrak{pp}$ of $\mathscr{O}_{K_1}$ which are non-dyadic and such that $v_\mathfrak{pp}(\alpha_1)=0$. Note that $T$ is infinite. For each $\mathfrak{p}\in T$, pick $\lambda_\mathfrak{p}\in \mathfrak{p}\setminus\mathfrak{p}^2$. By construction, we have \[ v_\mathfrak{q}(\lambda_\mathfrak{p})=\Bigl\lbrace\begin{array}{ccc}1 & \mbox{ if }&\mathfrak{q}=\mathfrak{p}\cr 0 & \mbox{ if } & \mathfrak{q}\neq\mathfrak{p} \end{array}. \]
Now, reasoning as in the previous case shows that  for any pair of distinct prime ideals $\mathfrak{p},\mathfrak{q}\in T$, $\lambda_\mathfrak{p}$ and $\lambda_\mathfrak{q}$ have distinct classes in the quotient group, hence the result.
\end{proof}

\begin{rem}Let us keep the notation of the proof.

If $L$ is a field, the construction of the set of prime ideals  $T$ and of $(u_\mathfrak{p})_{\mathfrak{p}}$ may be simplified when $\mathscr{O}_K$ is a PID. In this case, one may simply take for $T$ to be $S$ (or an infinite subset of $S$), and for $u_\mathfrak{p}$ a generator of one of the prime ideals of $\mathscr{O}_K$ lying above $\mathfrak{p}$.
 \end{rem}

\begin{ex}
Take $F=\mathbb{Q}, a_1=0$ and $a_1=1$, so that $\alpha_1=-4, \alpha_2=-3$ and  $K=\mathbb{Q}[\sqrt{12}]=\mathbb{Q}[\sqrt{3}]$. Thus, $\mathscr{O}_K=\mathbb{Z}[\sqrt{3}]$ is a PID.

If $p$ is a prime number such that $p\equiv -1 \ [12]$, then $p\equiv 3 \ [4]$, so $-1$ is not a square modulo $p$ and $p$ stays inert in $\mathbb{Q}[\sqrt{\alpha_1}]=\mathbb{Q}[i]$.

Now, we have \[ \genfrac(){}{0}{3}{p}=(-1)^{\frac{p-1}{2}\frac{3-1}{2}}\genfrac(){}{0}{p}{3}=
-\genfrac(){}{0}{p}{3}=-\genfrac(){}{0}{-1}{3}=1, \] so $p$ splits in $K$. 

Let $\pi_p\in \mathscr{O}_K$ be a generator of one of the prime ideals of $\mathscr{O}_K$ lying above $p$ (that is, we take $\pi_p\in\mathscr{O}_K$ such that $\vert N_{K/\mathbb{Q}}(\pi_p)\vert = p)$. The proof of the previous proposition and the remark above then show that the similarity classes of the bilinear forms corresponding to the classes of $\pi_p$, where $p$ runs  through the set of prime numbers satisfying $p\equiv - 1 \ [12]$, are pairwise distinct. Moreover, they are all indecomposable, since the class of $\pi_p$ is not trivial and not equal to the class of $\sqrt{3}$ (which can be seen by checking that $(-1,\pm 3p)$ is not split).

Using the explicit correspondence, one may show that the similarity class corresponding to $\pi_p=x+y\sqrt{3}=x+\dfrac{y}{2}\sqrt{12}$ is the class of the bilinear form $b_p$ on $\mathbb{Q}^2\otimes_\mathbb{Q} \mathbb{Q}^2$ such that   \[ \mathrm{Mar}(b_p; \mathscr{B}_0)=\begin{pmatrix}
    2x-2y& 2x+2y& 2x& 2x\cr -4y& 2x-2y& 0& 2x\cr  -2y& 2y& x-y& x+y \cr -4y& -2y& -2y& x-y
\end{pmatrix}. \] 
 \end{ex}

 The previous theorem shows that the determinant is not enough to distinguish decomposable bilinear spaces of dimension four with prescribed decomposable asymmetry, answering Question 2 negatively. However, the next theorem shows that the answer is positive at the level of anti-automorphisms. In particular, Question 1 has a positive answer, while Question 3 has a negative one.

 \begin{thm}\label{thmadjoint}
Let $(V',b') $ be a non-degenerate bilinear space of dimension four. Then $(\mathscr{L}(V') ,\sigma_{b'})$ is decomposable if and only if $b'$ has a decomposable asymmetry and trivial determinant.	
 \end{thm}
 
 \begin{proof}
 	Proposition \ref{adjclosed} shows that the two conditions above are necessary, so it remains to prove that they are also sufficient.
 	The case where $(V',b')$ has a non-generic decomposable asymmetry directly follows from Proposition \ref{non-gen}, so assume that the asymmetry of $(V',b')$ is conjugate to the asymmetry of $(V,b)=(F^2,b_{a_1})\otimes_F (F^2,b_{a_2})$, where $a_1,a_2\in F\setminus\{-2,2\}$ are such that $\alpha_i=a_i^2-4$ is not a square, for $i=1,2$.
 	
 	 In the sequel, we will make a slight abuse of notation, by using the same letter to denote an element of $\mathrm{Sym}(\mathscr{C}_b,*)^\times$ or its image in $K^\times$ via the isomorphism $\mathrm{Sym}(\mathscr{C}_b,*)^\times\simeq K^\times$.
 	  
    By Lemma \ref{lemdet}, $\det(b')$ is automatically trivial in this case, and we have to prove that $(\mathscr{L}(V'),\sigma_{b'})$ is indeed decomposable. By Proposition \ref{corres4}, $b'$ is similar to $b_u$ for some $u\in \mathrm{Sym}(\mathscr{C}_b,*)^\times$.

 	Let $u=r+s\sqrt{\alpha_1}\otimes \sqrt{\alpha_2}\in K^\times$.
 	If $s=0$, then the class of $u$ in $K^\times/F^\times N_{L/K}(L^\times)$ is trivial, and $b'$ is similar to $b$. Thus, $(V',b')$, and hence $(\mathscr{L}(V'),\sigma_{b'})$, are decomposable.
 	
 	We now assume that $s\neq 0$ for the rest of the proof.
 
 Let \[ E_1=I_2\otimes \begin{pmatrix}
 	0& \alpha_1 \cr 1 & 0
 \end{pmatrix}, \  F_1=\begin{pmatrix}
 0 & N_{K/F}(u)\cr 1 & 0 
 \end{pmatrix}\otimes \begin{pmatrix}
 1 & \phantom{-}0 \cr 0 & -1
 \end{pmatrix},\] \[ E_2= -\begin{pmatrix}
 1 & \phantom{-}0 \cr 0 & -1
 \end{pmatrix}\otimes \begin{pmatrix}
 0& \alpha_1 \cr 1 & 0
 \end{pmatrix}, \ F_2=-\begin{pmatrix}
 0 & N_{K/F}(u)\cr 1 & 0 
 \end{pmatrix}\otimes I_2.\]
 It is easy to see that the matrices $E_1,E_2,F_1,F_2\in \mathrm{M}_4(F)$ satisfy the following relations: \[ E_1^2=E_2^2=\alpha_1 I_4, \ F_1^2=F_2^2 =N_{K/F}(u)I_4, \ E_1F_1=-F_1E_1, \ E_2F_2=-F_2E_2, \]
 \[ E_1E_2=E_2E_1, \ E_1F_2=F_2E_1, \ F_1E_2=E_2F_1, \ F_1F_2=F_2F_1 \]
 
 It follows that there is a unique $F$-algebra isomorphism \[\Psi:(\alpha_1,N_{K/F}(u))\otimes_F(\alpha_1,N_{K/F}(u))\overset{\sim}{\to}\mathrm{M}_4(F)\] such that \[\Psi(i\otimes 1)=E_1, \ \Psi(j\otimes 1)=F_1, \ \Psi(1\otimes i)=E_2, \ \Psi(1\otimes j)=F_2.\]
 
 If $\gamma$ is the canonical involution on $(\alpha_1,N_{K/F}(u))$, one may check on the basis elements that we have $\Psi(\gamma\otimes \gamma)\Psi^{-1}=\mathrm{Int}(\Delta)^{-1}\circ^t$, where \[\Delta=\begin{pmatrix}
 	1 & 0 & 0 & 0\cr 0 & -\alpha_1 & 0 & 0 \cr 0 & 0 & -N_{K/F}(u) & 0 \cr 0 & 0 & 0 & \alpha_1 N_{K/F}(u)
 \end{pmatrix}. \]

Let $\sigma_1=\mathrm{Int}(a_1-2+i)\circ\gamma$ and $\sigma_2=\mathrm{Int}(a_2-2+\frac{1}{s\alpha_1}(ri+ij))\circ \gamma$. Finally, set $C=\Delta\big(\Psi((a_1-2+i)\otimes (a_2-2+\frac{1}{s\alpha_1}(ri+ij)))\big)^{-1}$. By construction, we have $\Psi(\sigma_1\otimes \sigma_2)\Psi^{-1}=\mathrm{Int}(C)^{-1}\circ^t$.

In particular, $\Psi$ induces an isomorphism of $F$-algebras with anti-automorphisms \[ ((\alpha_1,N_{K/F}(u)),\sigma_1)\otimes_F((\alpha_1,N_{K/F}(u)),\sigma_2)\simeq (\mathrm{M}_4(F), \mathrm{Int}(C)^{-1}\circ^t).  \]
We are going to show that the bilinear map $c$  on $F^2\otimes_F F^2$, whose matrix representation with respect to $\mathscr{B}_0$ is $C$, is similar to $b_u$, hence to $b'$, proving that  $(\mathscr{L}(V'),\sigma_{b'})$ is decomposable.

For, it is enough to show that the similarity class of $c$ corresponds to the class of $u$. One may check, using a computer algebra system, that we have the equality $C^{-1}C^t=P(A_{a_1}\otimes A_{a_2})P^{-1}$, where \[P=	\frac{1}{4s\alpha_1\alpha_2}\begin{pmatrix}
4s\alpha_1\alpha_2 & 2s\alpha_1\alpha_2 & 2s\alpha_1\alpha_2 & -r(a_1+2)(a_2+2)+s\alpha_1\alpha_2\cr
0 & 2r(a_2+2)& -2s\alpha_2(a_1+2) & r(a_2+2)-s\alpha_2(a_1+2)\cr 0 & 0 & 0 & -(a_1+2)(a_2+2)\cr 
0 &  2(a_2+2)& 0 & a_2+2
\end{pmatrix}.\]
Therefore, $c$ and $b$ have conjugate asymmetries, and the matrix representation $U_c$ with respect to $\mathscr{B}_0$ of a corresponding symmetric element $u_c\in \mathrm{Sym}(\mathscr{C}_b,*)^\times$ is $U_c=(P^tCP)^{-1} B$. One may check that we have 
\[U_c=-\frac{16s(a_1-2)(a_2-2)}{N_{K/F}(u)}(rI_4+ s (-a_1I_2+2A_{a_1})\otimes (-a_2I_2+2A_{a_2})).\]
After the appropriate identification, we see that this corresponds to the class of $u$ in $K^\times/F^\times N_{L/K}(L^\times)$. This concludes the proof.
 \end{proof}

\begin{rem}\label{q1q2}
The results of this section show that Question 1 and Question 2 are not equivalent, but let us emphasize this by showing what happens at the level of anti-automorphisms.

 Proposition \ref{critriv} shows that, if the class of $u$ is not trivial in $K^\times/F^\times N_{L/K}(L^\times)$, the quaternion $F$-algebra $(\alpha_1,N_{K/F}(u))$ is not split, hence is a division algebra. Now, if the class of $u$ is not equal to the trivial class or to the class of $\sqrt{\alpha_1}\otimes\sqrt{\alpha_2}$, then $(F^2\otimes_F F^2,b_u)$ is not decomposable by Proposition \ref{propdecomp}. Hence, in this case, $(\mathscr{L}(F^2\otimes_F F^2),\sigma_{b_u})$ decomposes as the tensor product of two division quaternion $F$-algebras with anti-automorphisms, but cannot be decomposed as the tensor product of two split quaternion $F$-algebras with anti-automorphisms, since $b_u$ would be decomposable otherwise. 
 \end{rem}

 Despite the fact that the determinant is not enough to detect decomposability of bilinear spaces, one may wonder if there are cohomological invariants of higher degrees which may be suitable for this purpose.

 We now proceed to show that it is not the case in general. First, we need to reinterpret our objects in a cohomological context.

\section{Galois Cohomology of bilinear spaces}\label{sec-coh}

We return for a moment to the general case. Let $F$ be a field of arbitrary characteristic, and let $(V,b)$ be a non-degenerate $F$-bilinear space.

In the following, the words ``algebraic group over $F$'' means ``affine algebraic group scheme over $F$'' in the sense of \cite{Wa}, that is a functor from the  category of commutative (associative unital) $F$-algebras to the category of groups, which may be represented by a finitely generated commutative $F$-algebra.

We denote by $\textbf{O}(b)$ and $\textbf{GO}(b)$ the algebraic group of automorphisms and similitudes of $(V,b)$ respectively.

We now use the techniques introduced in \cite{Ba} and pursued in \cite{Co} to introduce the point of view of algebras with involutions.

Keeping the notation of Section \ref{sec-corres}, it is easy to check that for all $f\in \mathscr{L}(V)$, we have $f\in \textbf{O}(b)(F)$ if and only if $\sigma_b(f)f=\Id_V$.
Now, if $f \in\textbf{O}(b)(F)$, then $f$ commutes with $a_b$. Indeed, for all $x,y,\in V$, we have \[  b(y,a_b(f(x)))=b(f(x),y)=b(x,f^{-1}(y))=b(f^{-1}(y),a_b(x))=b(y,f(a_b(x))), \]
and non-degeneracy implies the desired conclusion.

Therefore, $f\in \textbf{O}(b)(F)$ if and only if $f\in \textbf{U}(\mathscr{C}_b,*)(F)$.
The same is true if we replace $F$ by a commutative $F$-algebra $R$,
and $a_{b_R}=(a_b)_R$, and $\sigma_{b_R}=(\sigma_b)_R$.

We then have an equality of algebraic groups \[  \textbf{O}(b)=\textbf{U}(\mathscr{C}_b,*). \]
Similarly, we have  \[  \textbf{GO}(b)=\textbf{GU}(\mathscr{C}_b,*). \]

\begin{rem}
The idea to associate a certain algebra with involution to a bilinear form is not new. It has been introduced in \cite{Ba} in the context of systems of quadratic forms to reinterpret the automorphism group of such a system as a unitary group, and 
the equality between the group of similitudes of a bilinear form $b$ and  the group of similitudes of $(\mathscr{C}_b, *)$ already appears in \cite{Co}.

%Contrary to what is done in the papers \cite{Ba} and \cite{Co}, we will make only a limited use of this point of view. 
%It will only be needed to reinterpret Theorem \ref{corres} in cohomological terms in the generic case, and to establish the %smoothness of $\textbf{O}(b)$ and $\textbf{GO}(b)$ in characteristic different from two, using the well-known structure of %unitary groups (see the next lemma).
%
%Besides, it would be nice to have a direct proof of this smoothness property.
 \end{rem}

\begin{lem}\label{smooth}
Let $F$ be a field of characteristic different from two. Then $\mathbf{O}(b)$ and $\mathbf{GO}(b)$ are smooth.
\end{lem}

\begin{proof}
The proof of \cite[Lemma 6]{Wa} may be carried at the level of algebraic groups. Hence, by loc.cit., applied to $\textbf{U}(\mathscr{C}_b,*)$, we have an exact sequence of algebraic groups \[ 1\to \mathbf{N}\to \mathbf{O}(b)\to \mathbf{U'}\to 1, \]
where $\mathbf{N}$ is a split unipotent algebraic group, and $\mathbf{U'}$ is a direct product of Weil restrictions of classical groups. In particular, $\mathbf{N}$ and $\mathbf{U'}$ are smooth, and thus so is $\mathbf{O}(b)$.

Now, for every commutative $F$-algebra $R$, we have by definition an exact sequence \[ 1\to \mathbf{O}(b)(R)\to \mathbf{GO}(b)(R)\to R^\times. \]
Moreover, for any $\lambda\in \overline{F}^\times$, $\sqrt{\lambda}\Id_V$ is a similitude of $b_{\overline{F}}$ with similarity factor equal to $\lambda$. Hence the map $\mathbf{GO}(b)(\overline{F})\to \overline{F}^\times$ is surjective. Since $\mathbb{G}_{m,F}$ is smooth, 
by \cite[Proposition (22.10)]{KMRT}, we get an exact sequence of algebraic groups \[ 1\to \mathbf{O}(b)\to \mathbf{GO}(b)\to \mathbb{G}_{m,F}\to 1 \]
The smoothness of $\mathbf{GO}(b)$ then comes from the smoothness of $\mathbf{O}(b)$ and $\mathbb{G}_{m,F}$.
\end{proof}

We may now prove the following proposition.

\begin{prop}\label{cohogo}
Let $F$ be a field of characteristic different from two, and let $(V,b)$ be a non-degenerate bilinear space.

Then the pointed set $H^1(F, \mathbf{O}(b))$, resp. $H^1(F, \mathbf{GO}(b))$, is in one-to-one correspondence with the pointed set of isomorphism classes, resp. of similarity classes, of bilinear spaces $(V',b')$ such that $a_{b'}$ and $a_b$ are conjugate (the base point being the isomorphism class, resp. the similarity class, of $(V,b)$).
\end{prop}

\begin{proof}
The groups $\mathbf{O}(b)$ and $\mathbf{GO}(b)$ may be viewed as the stabilizer of $(V,b)$ under the obvious actions of 
$\textbf{GL}(V)$ and $\textbf{GL}(V)\times\mathbb{G}_m$ on bilinear spaces respectively. Since  $\textbf{GL}(V)$ and $\textbf{GL}(V)\times\mathbb{G}_m$ have trivial $H^1$,
Galois descent shows that the pointed set $H^1(F, \mathbf{O}(b))$, resp. $H^1(F, \mathbf{GO}(b))$, is in one-to-one correspondence with the pointed set of isomorphism classes, resp. of similarity classes, of bilinear spaces $(V',b')$ which is isomorphic, resp. similar, to $(V,b)$ over $F_{sep}$.
Now, since $\mathbf{O}(b)$ and $\mathbf{GO}(b)$ are smooth by Lemma \ref{smooth},  this last condition is equivalent to ask for  $(V',b')$ to be isomorphic, resp. similar, to $(V,b)$ over $\overline{F}$ (see \cite{Wat}, Section 18.5).
By Lemma \ref{symfbar}, we are done.
\end{proof}

\begin{rem}
In view of the previous result, the first part of Theorem \ref{corres} may then be recovered when $F$ has characteristic different from two when $\mathscr{C}_b$ is separable using the results of Section 29.D of \cite{KMRT} (these results are proved in the context of central simple algebras with involution, but may be easily extended to the separable case).
 \end{rem}

We now give a cohomological reinterpretation of some results of the previous section.

Assume until the end of this section that $F$ has characteristic different from two, and that $b=b_{a_1}\otimes b_{a_2}$, where $a_1,a_2\in F\setminus\{2,-2\}$ are such that $\alpha_i=a_i^2-4$ is not a square, for $i=1,2$.

Recall that  $L=F[\sqrt{\alpha_1}]\otimes_F F[\sqrt{\alpha_2}]$ and that $K=F[\sqrt{\alpha_1}\otimes \sqrt{\alpha_2}]$.

Hence, we have the following isomorphisms of algebraic group schemes
\[ \textbf{O}(b)=\textbf{U}(\mathscr{C}_b,*)\simeq \textbf{U}(L,\sigma)\simeq R_{K/F}(\mathbb{G}^{(1)}_{m,L/K}), \]
where $\mathbb{G}^{(1)}_{m,L/K}=\ker(R_{L/K}(\mathbb{G}_{m,L})\overset{N_{L/K}}{\longrightarrow} \mathbb{G}_{m,K})$.

We then recover the smoothness of $\textbf{O}(b)$, and thus the smoothness of $\textbf{GO}(b)$, in this particular case. We may also recover the result of Theorem \ref{corres} by applying Galois cohomology, as follows. 

Since $\textbf{GO}(b)=\textbf{GU}(\mathscr{C}_b,*)\simeq \textbf{GU}(L,\sigma)$, we have an exact sequence 
\[ 1\to \textbf{GO}(b)\to R_{L/F}(\mathbb{G}_{m,L})\times \mathbb{G}_{m,F}\to R_{K/F}(\mathbb{G}_{m,K})\to 1,  \]
where the second map is defined on the $R$-points by \[ (x,\lambda)\in (L\otimes_FR)^\times\times R^\times\mapsto \lambda^{-1} N_{L\otimes_F R/K\otimes_F R}(x). \]

Applying Galois cohomology, Shapiro's lemma and Hilbert 90 yield a group exact sequence \[ L^\times \times F^\times \to K^\times\to H^1(F,\textbf{GO}(b))\to 1, \]
and the first isomorphism theorem then gives \[ H^1(F, \textbf{GO}(b))\simeq K^\times/F^\times N_{L/K}(L^\times). \]

\section{Cohomological invariants and decomposability}\label{sec-cohinv}
 
We now address the question of detecting decomposability of bilinear spaces in terms of the vanishing of some cohomological invariants.

Let $\mathscr{H}$ be a functor from the category of field extensions of $F$ to the category of abelian groups. If $G$ is an algebraic group over $F$, a {\it normalized invariant of $G$ with values in} $\mathscr{H}$ is a natural transformation $\eta: H^1(_-, G)\to\mathscr{H}$ of functors from the category of field extensions of $F$ to the category of pointed sets.

When  $C$ is a discrete $\mathrm{Gal}(F_{sep}/F)$-module and $\mathscr{H}=H^d(_-,C)$, such an invariant is called a \textit{normalized cohomological invariant of $G$ of degree $d$ with coefficients in $C$}.

We denote by $\mathrm{Inv}_0^d(G,C)$ the group of normalized cohomological invariants of $G$ of degree $d$ with coefficients in $C$.

We have already seen that the determinant detects decomposability of bilinear spaces with a non-generic asymmetry, so we will focus on the generic case.

In the sequel, $(V,b)$ will be the $F$-bilinear space $(F^2,b_{a_1})\otimes_F (F^2,b_{a_2})$, where $a_1,a_2\in\setminus\{-2,2\}$ are such that $\alpha_i=a_i^2-4$ is not a square, for $i=1,2$.

As before, we set $L=F[\sqrt{\alpha_1}]\otimes_F F[\sqrt{\alpha_2}]$, and $K=F[\sqrt{\alpha_1}\otimes\sqrt{\alpha_2}]$.

If $E/F$ is a field extension, we will denote $E\otimes_F K$ and $E\otimes_F L$ by $K_E$ and $L_E$ respectively.

From now on, if $u\in K_E^\times$, we will denote by $[u]$ its class in $K_E^\times/E^\times N_{L_E/K_E}(L_E^\times)$. In view of the results of Section \ref{sec-det-gen}, we will identify $H^1(E,\textbf{GO}(b))$ with this group, and we will often write $[u]\in H^1(E, \textbf{GO}(b))$.

Let us start by giving some examples of invariants of $\mathbf{GO}(b)$.

\begin{ex}\label{exnu}
If $E/F$ is a field extension, let $\nu_E: H^1(E, \mathbf{GO}(b))\to H^2(E,\mathbb{Z}/2\mathbb{Z})$ be the map defined by \[  \nu_E([u])=
(\alpha_1)\cup(N_{K_E/E}(u))  \mbox{ for all }u\in K_E^\times.\]

Proposition \ref{critriv}, together with the usual identification $H^2(E)\simeq \mathrm{Br}_2(E)$, shows in particular that $\nu_E$ is well-defined.

The various maps $\nu_E$ are compatible with scalar extension and thus define a normalized cohomological invariant $\nu$ of $\mathbf{GO}(b)$ of degree 2 with coefficients in $\mathbb{Z}/2\mathbb{Z}$.

More generally, if $d\geq 2$ and $\rho\in H^{d-2}(F,\mathbb{Z}/2\mathbb{Z})$, one may define a normalized cohomological invariant $\rho\cup \nu$ of degree $d$ with coefficients in $\mathbb{Z}/2\mathbb{Z}$ by setting  \[ (\rho\cup\nu)_E([u])=\rho_E\cup\nu_E([u])=\rho_E\cup (\alpha_1)\cup(N_{K_E/E}(u)) \ \mbox{ for all }u\in K_E^\times.\]
 \end{ex}

\begin{rem}\label{remnu}
	Assume that $\alpha_2\alpha_1^{-1}$ is a square. By Remark \ref{a1a2square}, for any field extension $E/F$, we have a group isomorphism 
	\[K_E^\times/E^\times N_{L_E/K_E}(L_E^\times)\simeq E^\times/N_{E\otimes_F K_1/E}((E\otimes_F K_1)^\times),\] where $K_1=F[\sqrt{\alpha_1}]$.
	The last part of this remark show that, using this isomorphism, the map $\nu_E$ sends the class of $\lambda $ in $ E^\times/N_{E\otimes_F K_1/E}((E\otimes_F K_1)^\times)$ to $(\alpha_1)\cup (\lambda)$.
 \end{rem}

\begin{ex}\label{nutilde}
Let $\nu$ be the normalized cohomological invariant of degree two defined in Example \ref{exnu}. 
For any field extension $E/F$, composing with the canonical projection yields a map 	\[\widetilde{\nu}_E:H^1(E, \mathbf{GO}(b))\to H^2(E,\mathbb{Z}/2\mathbb{Z})/\langle (\alpha_1)_E\cup (\alpha_2)_E\rangle. \]

These maps fit together to define an invariant (with obvious notation) \[\widetilde{\nu}:H^1(_-,\mathbf{GO}(b))\to H^2(-,\mathbb{Z}/2\mathbb{Z})/\langle (\alpha_1)\cup (\alpha_2)\rangle. \] 
 \end{ex}

 Proposition \ref{critriv} shows immediately that, for every field extension, the map $\nu_E$ is injective. Since $\nu_E(b_E)=0$ and $\nu_E((b_{-a_1}\otimes b_{-a_2})_E)=(\alpha_1)_E\cup(\alpha_2)_E$, it easily follows that the map \[\widetilde{\nu}_E:H^1(E, \mathrm{GO}(b))\to H^2(E)/\langle (\alpha_1)_E\cup (\alpha_2)_E\rangle\] vanishes exactly at $b_E$ and $(b_{-a_1}\otimes b_{-a_2})_E$.
Consequently, we have the following proposition.

\begin{prop}\label{dec-detect}
The invariant $\widetilde{\nu}:H^1(_-,\mathbf{GO}(b))\to H^2(-,\mathbb{Z}/2\mathbb{Z})/\langle (\alpha_1)\cup (\alpha_2)\rangle $ detects decomposability.

In other words, for any field extension $E/F$, and any non-degenerate $(V',b')$ $E$-bilinear space whose asymmetry is similar to the asymmetry of $(V,b)_E$, we have $\nu_E(b')=0$ if and only if $(V',b')$ is decomposable.
\end{prop}	

\begin{rem}
Even if $\widetilde{\nu}$ is not a cohomological invariant per se, we might consider that Proposition \ref{dec-detect} settles the problem of detecting decomposability in the case of bilinear spaces with a generic decomposable asymmetry. However, the construction of $\widetilde{\nu}$ is a bit too ad-hoc, and the class $(\alpha_1)\cup(\alpha_2)$ seems to appear from nowhere.
It would be nice to have a conceptual construction of $\widetilde{\nu}$, which could explain why we have to mod out by the class $(\alpha_1)\cup(\alpha_2)$.
 \end{rem}
	
In view of the previous remark, one may wonder if we could find a family of cohomological invariants with coefficients in a discrete Galois module $C$ which detect decomposability. We start with a lemma.

\begin{lem}\label{tor2}
Any normalized cohomological invariant of $\mathbf{GO}(b)$ of degree $d$ with coefficients in $C$ is killed by two. 
\end{lem} 

\begin{proof}
Let $\eta\in \mathrm{Inv}_0^d(\textbf{GO}(b),C)$, let $E/F$ be a field extension, and let $u\in K_E^\times$. 

If $\alpha_1$ is a square in $E$, then $L_E$ is the split quadratic \'{e}tale $K_E$-algebra, and the norm map $N_{L_E/K_E}$ is surjective, so $u\sim 1$. In this case, we obviously have $2\cdot\eta_E([u])=0$.

Assume that $\alpha_1$ is not a square in $E$, so that $E'=E\otimes_F F[\sqrt{\alpha_1}]$ is a field.
Then, $L_{E'}$ is the split quadratic \'{e}tale $K_{E'}$-algebra, so the norm map  $N_{L_{E'} /K_{E'}}$ is surjective. 
It follows  that the image $[u]_{E'}$ of $[u]$ under the restriction map \[ H^1(E, \textbf{GO}(b))\to H^1(E', \textbf{GO}(b)) \] is trivial. Since $\eta$ is normalized, it yields \[ (\eta_E([u]))_{E'}=\eta_{E'}([u]_{E'})=\eta_{E'}([1])=0. \]
 Since $[E':E]=2$, composing with the corestriction map shows again that 
 $2\cdot \eta_E([u])=0$. 
\end{proof}

The previous lemma shows that we have to consider Galois modules $C$ having $2$-torsion if we want to have a chance to 
get non-trivial cohomological invariants.

We will stick to the case $C=\mathbb{Z}/2\mathbb{Z}$, and  write \[ H^d(F)  \mbox{ and }  \mathrm{Inv}_0^d(\textbf{GO}(b)) \] instead of  
$H^d(F,\mathbb{Z}/2\mathbb{Z})$ and $\mathrm{Inv}_0^d(\textbf{GO}(b),\mathbb{Z}/2\mathbb{Z})$.

Since $\textbf{GO}(b)$ is a connected algebraic group in our case, any element $\mathrm{Inv}_0^1(\textbf{GO}(b))$
is necessarily trivial by \cite[Proposition 31.15]{KMRT}, so we may assume that $d\geq 2$.

Note that since the determinant may be interpreted as an element of $\mathrm{Inv}_0^1(\textbf{GO}(b))$, we recover Lemma \ref{lemdet}.

The next theorem shows that all the normalized cohomological of $\textbf{GO}(b)$ of degree $d$ have the form $\rho \cup \nu$, as defined in Example \ref{exnu}, under a mild assumption on $F$ if $d\geq 3$ and $L$ is a field.

\begin{thm}\label{thmcoh}
Let $F$ be a field of characteristic different from two, and let $d\geq 2$. 
If $d\geq 3$ and $L$ is a field, assume that every element of $F$ is a sum of two squares.

Then, for all $\eta\in \mathrm{Inv}_0^d(\textbf{GO}(b))$, there exists $\rho\in H^{d-2}(F)$ such that $\eta=\rho\cup \nu$.	
	\end{thm}

	\begin{proof}
	Assume first that $\alpha_2\alpha_1^{-1}$ is a square. For every field extension $E/F$, $H^1(E,\mathbf{GO}(b))$ then identifies to 	 $E^\times/N_{E\otimes_F K_1/E}((E\otimes_F K_1)^\times)$ by Remark \ref{a1a2square}. It follows that the class $[t]\in F(t)^\times/N_{K_1(t)/F(t)}(K_1(t)^\times)$
	is a generic $\mathbf{GO}(b)$-torsor.	
	
	Let  $\eta\in \mathrm{Inv}_0^d(\textbf{GO}(b))$, where $d\geq 2$. Since $[t]$ is unramified outside $t$, so is $\eta_{F(t)}([t])$ 
	 by \cite[Thm 11.7]{GMS}. It follows that  there exist $\rho_0\in H^d(F)$ and  $\rho_1\in H^{d-1}(F)$ such that  $\eta_{F(t)}([t])=(\rho_0)_{F(t)}+(\rho_1)_{F(t)}\cup (t)$ by  \cite[§ 9.4]{GMS}. Since $\eta$ is normalized, evaluation at $t=0$ yields $\rho_0=0$.
	 Hence, $\eta_{F(t)}([t])=(\rho_1)_{F(t)}\cup (t)$.
	 Now, since $K_1\otimes_F K_1(t)$ is isomorphic to $K_1(t)\times K_1(t)$ as a $K_1(t)$-algebra, the norm map $N_{K_1\otimes_F K_1(t)}/K_1(t)$ is surjective, and $[t]_{K_1(t)}$ is the trivial class. Since $\eta$ is normalized, it follows that we have 
	 \[ (\eta_{F(t)}([t]))_{K_1(t)}=\eta_{K_1(t)}([1])=0=(\rho_1)_{K_1(t)}\cup(t). \]
	 Evaluation at $t=0$ yields $(\rho_1)_{K_1}=0$. Therefore, there exists $\rho\in H^{d-2}(F)$ such that $\rho_1=\rho\cup (\alpha_1)$.
	 All in all,  $\eta_{F(t)}([t])=\rho_{F(t)}\cup (\alpha_1)\cup (t)=(\rho\cup \nu)_{F(t)}([t])$.
	Thus, the cohomological invariants $\eta$ and $\rho\cup \nu$ coincide on a versal torsor, hence are equal by \cite[Thm 12.3]{GMS} (this uses the smoothness of $\textbf{GO}(b)$).

	We now assume until the end of the proof that $L$ is a field. We will then write $\sqrt{\alpha_1\alpha_2}$ instead of $\sqrt{\alpha_1}\otimes \sqrt{\alpha_2}$.

We start by checking that the class $[1+t\sqrt{\alpha_1\alpha_2}]\in H^1(F(t),\textbf{GO}(b))$ is a versal $\textbf{GO}(b)$-torsor.

We have to show that, for any field extension $E/F$, and any $u\in K_E^\times$, there exists $s\in E$ such that $u\sim 1+ s\sqrt{\alpha_1\alpha_2}$.
Abusing notation, we will write $1$ and $\sqrt{\alpha_1\alpha_2}$ instead of $1\otimes 1$ and $1\otimes\sqrt{\alpha_1\alpha_2}\in K_E$.

If $u=x+y\sqrt{\alpha_1\alpha_2}$ with $x\in E^\times$, one may simply take $s=yx^{-1}$.
If $x=0$, we have 
\[ u\sim \sqrt{\alpha_1\alpha_2}\sim \dfrac{\sqrt{\alpha_1\alpha_2}}{2\alpha_1\alpha_2}N_{L_E/K_E}(1+\sqrt{\alpha_1\alpha_2}), \] that is \[ u\sim \dfrac{\sqrt{\alpha_1\alpha_2}(1+\sqrt{\alpha_1\alpha_2})^2}{2\alpha_1\alpha_2}=1+\dfrac{1+\alpha_1\alpha_2}{2\alpha_1\alpha_2}\sqrt{\alpha_1\alpha_2}. \]

This yields the desired conclusion.

We are now ready to determine the cohomological invariants of $\mathbf{GO}(b)$.
We start with the case $d=2$.

Let $\eta\in \mathrm{Inv}_0^2(\textbf{GO}(b))$ and set $\xi=\eta_{F(t)}([1+t\sqrt{\alpha_1\alpha_2}])$. 
We saw  during the proof of Lemma \ref{tor2} that $\xi_{F(t)[\sqrt{\alpha_1}]}=0$, so we get $\xi=(\alpha_1)\cup \beta$, where $\beta\in H^1(F(t)). $
Moreover, $[1+t\sqrt{\alpha_1\alpha_2}]$ is unramified outside $\pi=1-\alpha_1\alpha_2 t^2$, and thus, so is $\xi$ by \cite[Thm 11.7]{GMS}.
We will identify the residue field corresponding to the $\pi$-adic valuation to $K$, by sending the class of $t$ modulo $\pi$ onto $\dfrac{1}{\sqrt{\alpha_1\alpha_2}}$. 

 We may write $\beta=(P \pi^k)$, for some $P\in F[t]\setminus\{0\}$ not divisible by $\pi$, and $k=0$ or $1$, so $\partial_\pi(\xi)=k(\alpha_1)\in H^1(K)$. 

Therefore, $\xi+k(\alpha_1)\cup(\pi)$ is unramified at $\pi$. It is also unramified outside $\pi$ because both summands are. Finally,  $\xi+k(\alpha_1)\cup(\pi)$ is unramified, hence constant by \cite[Thm 10.1]{GMS}. Thus, we get \[ \xi=\beta'_{k(t)}+k(\alpha_1)\cup(\pi)=\beta'_{k(t)}+k(\alpha_1)\cup(1-\alpha_1\alpha_2 t^2) \] for some $\beta'\in H^2(F)$.
The evaluation of $\xi$ at $t=0$ is trivial, since $\eta$ is a normalized cohomological invariant, thus we get $\beta'=0$ and \[ \xi=\eta_{F(t)}([1+t\sqrt{\alpha_1\alpha_2}])=k(\alpha_1)\cup(1-\alpha_1\alpha_2t^2)=k\nu_{F(t)}([1+t\sqrt{\alpha_1\alpha_2}]). \]

Thus, the cohomological invariants $\eta$ and $k \nu$ coincide on a versal torsor, hence are equal by \cite[Thm 12.3]{GMS}.

We now assume that $d\geq 3$, and that every element of $F$ is a sum of two squares.

Let $\eta\in \mathrm{Inv}_0^d(\textbf{GO}(b))$, where $d\geq 3$, and set  $\xi=\eta_{F(t)}([1+t\sqrt{\alpha_1\alpha_2}]).$
As in the previous case, we have  $\xi=(\alpha_1)\cup \beta$, where $\beta\in H^{d-1}(F(t))$,
and $\xi$ is unramified outside $\pi=1-\alpha_1\alpha_2 t^2$.
Note that $\partial_\pi(\xi)=(\alpha_1)\cup \partial_\pi\beta\in H^{d-1}(K).$ 
Since $K/F$ is a quadratic extension, the residue of $\beta$ at $\pi$ has the form  \[ \rho_K+\gamma_K\cup(-r+s\frac{1}{\sqrt{\alpha_1\alpha_2}}),\] where $(r,s)\neq (0,0)\in F\times F$, $\rho\in H^{d-2}(F)$ and $\gamma\in H^{d-3}(F)$. 

If $s=0$ (and thus $r\neq 0$), replacing $\rho$ by $\rho+\gamma\cup (-r)$, we may assume that $\gamma=0$.
In this case, $\xi+ \rho_{F(t)}\cup (\alpha_1)\cup (\pi)$ is unramified at $\pi$. It is also unramified outside $\pi$ since both summands are.
Hence, it unramified, hence constant. Taking into account that $\eta$ is normalized, specialization at $t=0$ shows that \[ \xi+ \rho_{F(t)}\cup (\alpha_1)\cup (\pi)=0, \] that is $\xi=\rho_{F(t)}\cup (\alpha_1)\cup (\pi)$. In other words, \[  \eta_{F(t)}([1+t\sqrt{\alpha_1\alpha_2}])=(\rho\cup\nu)_{F(t)}
([1+t\sqrt{\alpha_1\alpha_2}]). \]
Again, this implies the equality $\eta=\rho\cup \nu$.

We now assume until the end that $s\neq 0$. Replacing $\rho$ by $\rho+\gamma\cup (s)$ and $r$ by $\dfrac{r}{s}$, we may assume that $s=1$.
In this case, we have \[ \partial_\pi(\xi)=(\alpha_1)\cup \rho_K+(\alpha_1)\cup\gamma_K\cup (-r+\frac{1}{\sqrt{\alpha_1\alpha_2}}), \]
and
$\xi+ (\alpha_1)\cup \rho_{F(t)}\cup (\pi)+(\alpha_1)\cup \gamma_{F(t)}\cup (t-r )\cup (\pi)$ is unramified at $\pi$. It is also unramified outside $t-r $, since all summands are. Now the residue of this last cohomology class at $t-r$ is 
$(\alpha_1)\cup \gamma\cup(\pi(r))$.
All in all, the class \[ \xi+ (\alpha_1)\cup \rho_{F(t)}\cup (\pi)+(\alpha_1)\cup \gamma_{F(t)}\cup [ (t-r)\cup (\pi)+ (t-r)\cup  (\pi(r))] \] is unramified, hence constant. Therefore, there exists $\delta\in H^d(F)$ such that
\[ \xi=\delta_{F(t)}+ (\alpha_1)\cup \rho_{F(t)}\cup (\pi)+(\alpha_1)\cup \gamma_{F(t)}\cup [(t-r)\cup (\pi)+ (t-r)\cup  (\pi(r))], \]
Note that, if $r=0$, since $\pi(0)=1$, and thus $(\pi(0))=0$, the last term is $0$.
Now, the specialization of $\xi$ at $t=0$ is trivial, since $\eta$ is normalized, hence we get 
\[  0=\delta+ (\alpha_1)\cup \gamma\cup(-r)\cup(\pi(r)), \] and finally 
\[ \xi=(\alpha_1)\cup \rho_{F(t)}\cup (\pi)+(\alpha_1)\cup \gamma_{F(t)}\cup [(t-r)\cup (\pi)+ (-r(t-r))\cup  (\pi(r))], \]
where the last term has to be interpreted as $0$ if $r=0$.

By the properties of versal torsors, it follows that, for any field extension $E/F$, and any $y\in E$, we have 
\[ \eta_E([1+y\sqrt{\alpha_1\alpha_2}])=\begin{array}{lll}(\alpha_1)\cup \rho_E\cup (\pi(y))&+&(\alpha_1)\cup \gamma_E\cup [(y-r)\cup (\pi(y))\cr &+& (-r(y-r))\cup  (\pi(r))] \ \mbox{ (1) }.\end{array} \]

Let $X,Y$ be two indeterminates over $F$. Set $F'=F(X,Y)$ and \[ \xi'=\eta_{F'}([1+X\sqrt{\alpha_1\alpha_2}][1+Y\sqrt{\alpha_1\alpha_2}])\in H^d(F'). \]

Note that $\xi'$ is unramified outside $\pi(X)$ and $\pi(Y)$.
We have \[ (1+X\sqrt{\alpha_1\alpha_2})(1+Y\sqrt{\alpha_1\alpha_2})=1+\alpha_1\alpha_2 XY+(X+Y)\sqrt{\alpha_1\alpha_2}, \] hence 
$[(1+X\sqrt{\alpha_1\alpha_2})(1+Y\sqrt{\alpha_1\alpha_2})]=[ 1+\frac{X+Y}{1+\alpha_1\alpha_2 XY}\sqrt{\alpha_1\alpha_2}]$, 
and also \[ 1-\alpha_1\alpha_2 \frac{(X+Y)^2}{(1+\alpha_1\alpha_2XY)^2}=\frac{\pi(X)\pi(Y)}{(1+\alpha_1\alpha_2 XY)^2}. \]

Therefore, we have \[ \xi'=\begin{array}{c}(\alpha_1)\cup \rho _{F'}\cup (\pi(X)\pi(Y))+ (\alpha_1)\cup \gamma _{F'}\cup [(\frac{X+Y}{1+\alpha_1\alpha_2 XY}-r)\cup (\pi(X)\pi(Y))+ \cr \vspace{6pt}(-r(\frac{X+Y}{1+\alpha_1\alpha_2 XY}-r)\cup (\pi(r))]\end{array}. \]

Set $\xi_X=\eta _{F'}([1+X\sqrt{\alpha_1\alpha_2}])$ and $\xi_Y=\eta _{F'}([1+Y\sqrt{\alpha_1\alpha_2}])$. Note that $\pi_X$ is unramified outside $\pi(X)$, while $\pi_Y$ is unramified outside $\pi(Y)$.

By equality $(1)$, we have \[  \xi_X=(\alpha_1)\cup \rho _{F'}\cup (\pi(X))+(\alpha_1)\cup \gamma _{F'}\cup [(X-r)\cup (\pi(X))+ (-r(X-r))\cup (\pi(r))],  \] and 
\[  \xi_Y=(\alpha_1)\cup \rho _{F'}\cup (\pi(Y))+(\alpha_1)\cup \gamma _{F'}\cup [(Y-r)\cup (\pi(Y))+ (-r(Y-r))\cup (\pi(r))].  \]

We now proceed to show that $\xi'=\xi_X+\xi_Y$.
The cohomology class $\xi'+\xi_X+\xi_Y$ is unramified outside $\pi(X)$ and $\pi(Y)$.
Now, \[ \partial_{\pi(X)}(\xi')=(\alpha_1)\cup \rho_{K(Y)}+(\alpha_1)\cup\gamma_{K(Y)}\cup (\frac{\frac{1}{\sqrt{\alpha_1\alpha_2}}+Y}{1+Y\sqrt{\alpha_1\alpha_2}}-r),  \]
that is \[ \partial_{\pi(X)}(\xi')=(\alpha_1)\cup \rho_{K(Y)}+(\alpha_1)\cup\gamma_{K(Y)}\cup (\frac{1}{\sqrt{\alpha_1\alpha_2}}-r)\in H^{d-1}(K(Y)). \]
Moreover, $\partial_{\pi(X)}(\xi_X)=(\alpha_1)\cup \rho_{K(Y)}+(\alpha_1)\cup\gamma_{K(Y)}\cup (\frac{1}{\sqrt{\alpha_1\alpha_2}}-r)\in H^{d-1}(K(Y))$, while $\partial_{\pi(X)}(\xi_Y)=0$. It follows that  $\xi'+\xi_X+\xi_Y$ is unramified at $\pi(X)$. Exchanging the roles of $X$ and $Y$ shows that  $\xi'+\xi_X+\xi_Y$ is also unramified at $\pi(Y)$. Consequently,  $\xi'+\xi_X+\xi_Y$ is unramified, hence constant. Specialization at $X=0$ and $Y=0$ shows that  $\xi'+\xi_X+\xi_Y=0$, that is $\xi'=\xi_X+\xi_Y$.

In other words, we get \[ \eta _{F'}([1+X\sqrt{\alpha_1\alpha_2}][1+Y\sqrt{\alpha_1\alpha_2}])=\eta _{F'}([1+X\sqrt{\alpha_1\alpha_2}])+\eta _{F'}([1+Y\sqrt{\alpha_1\alpha_2}]). \]
A specialization argument then yields that $\eta_E$ is a group morphism for any field extension $E/F$.

Note that we have a surjection $\pi: R_{K/F}(\mathbb{G}_{m,K})\to H^1 (_-, \textbf{GO}(b))$. The map $\eta'=\eta\circ \pi$ is then a cohomological invariant of $R_{K/F}(\mathbb{G}_{m,K})$ is the sense of \cite{MPT}, that is, a natural transformation between functors from the category of field extensions of $F$ to the category of abelian groups.
By \cite[Thm 1.1]{MPT}, there exists $\varphi\in H^{d-1}(K)$ such that for any field extension $E/F$, and any $u\in K_E^\times$, we have 
\[ \eta'_E(u)=\mathrm{Cor}_{K_E/E}(\varphi_{K_E}\cup (u)). \]
Applying this to $E=F(t)$ and $u=t$ yields \[ \eta'_{F(t)}(t)=\eta_{F(t)}([t])=\mathrm{Cor}_{K(t)/F(t)}(\varphi_{K(t)})\cup (t)=0, \]
since $[t]=[1]$.
Now, $\varphi \in H^{d-1}(K)$, so it is unramified at $t$. Computing residues at $t$ shows that $\mathrm{Cor}_{K/F}(\varphi)=0$. Hence $\varphi=\psi_K$ for some $\psi\in H^{d-1}(F)$ and $\eta'_E(u)=\psi_E\cup (N_{K_E/E}(u))$ for any field extension $E/F$ and any $u\in K_E^\times$.

On the other hand, for all $\lambda \in L_E^\times$, we have \[ \eta'_E(N_{L_E/K_E}(\lambda))=\eta_E([N_{L_E/K_E}(\lambda)])=\eta_E([1])=0,  \]
that is \[ \psi_E\cup (N_{L_E /E}(\lambda))=0. \]

Note that $L(t)\otimes _F L\simeq L(t)^4$ as $L(t)$-algebras. In particular,  $N_{L(t)\otimes_F L/L(t)}$ is surjective.
Applying the previous equality to $E=L(t)$ then yields $\psi_{L(t)}\cup (\mu)=0$ for all $\mu\in L(t)^\times$. Taking $\mu=t$ and computing residues shows that $\psi_L=0$.

Since any element of $F$ is a sum of two squares, by \cite[Thm. 1]{Ka},  we have $\psi=(\alpha_1)\cup\psi_1+(\alpha_2)\cup \psi_2+ (\alpha_1\alpha_2)\cup\psi_3$, for some $\psi_i\in H^{d-2}(F)$.

Since $(\alpha_2)=(\alpha_1)+(\alpha_1\alpha_2)$, this equality may be rewritten as \[ \psi=(\alpha_1)\cup \rho+(\alpha_1\alpha_2)\cup \rho', \] for some $\rho,\rho'\in H^{d-2}(F)$.

We then finally get that $\eta'_E(u)=(\alpha_1)\cup\rho_E\cup (N_{K_E/E}(u))$ for any field extension $E/F$, and any $u\in K_E^\times$. that is, $\eta=\rho\cup \nu$.
\end{proof}

\begin{rem}
	The assumption on $F$ is necessary to apply Theorem 1 of \cite{Ka}.
	As pointed out by Bruno Kahn in his paper, this result is probably true unconditionally. Therefore, our theorem should be true unconditionally as well, that is,  for all $d\geq 2$, and for all $\eta\in \mathrm{Inv}_0^d(\textbf{GO}(b))$, there exists $\rho\in H^{d-2}(F)$ such that $\eta=\rho\cup \nu$.
 \end{rem}

 \begin{rem}\label{negquest2}
Theorem \ref{thmcoh} shows that characterizing decomposability using the vanishing of some cohomological invariants with coefficients in $\mathbb{Z}/2\mathbb{Z}$ is somewhat hopeless. 
 
 Indeed, assume that the quaternion $F$-algebra $(\alpha_1,\alpha_2)$ is not split. Then, Remark \ref{remdec} and Lemma \ref{decasym} show that 
 $b$ and $b'=b_{-a_1}\otimes b_{-a_2}$ are not similar, and are the only decomposable bilinear forms whose asymmetry is conjugate to $a_b$, up to similarity. However, if $H^3(F)=0$, there is no normalized cohomological invariant of $\textbf{GO}(b)$ with coefficients in $\mathbb{Z}/2\mathbb{Z}$ which vanishes exactly at $b$ and $b'$.

As a matter of fact, if $d=2$, Theorem \ref{thmcoh} shows that $\eta=0$ or $\nu$. But, we have \[  \nu_F(b')=(\alpha_1)\cup (\alpha_2)\neq 0 \]
by assumption.
If $H^3(F)=0$, then $H^d(F)=0$ for all $d\geq 3$ and thus $\eta_F=0$ for all $\eta\in \mathrm{Inv}^d_0(\textbf{GO}(b))
$. We then get the desired result.
In fact, the assumption $H^3(F)=0$ even implies that every cohomological invariant of degree $d\geq 5$ is identically zero, provided that every element of $F$ is a sum of two squares. This will be the case for example if $F$ is a number field containing~$i$. 
This answers Question 4 negatively.

However, one may wonder if one could get examples of a field $F$ such that, for all $d\geq 2$, there exists at least one $\eta\in \mathrm{Inv}_0^d(\textbf{GO}(b))$ such that $\eta_F\neq 0$, and for which there is still no normalized cohomological invariant of $\textbf{GO}(b)$ with coefficients in $\mathbb{Z}/2\mathbb{Z}$ which vanishes exactly at $b$ and $b'$. This is the goal of the next  examples.
  \end{rem}

\begin{ex}\label{laurentr}
	Let $k$ be a real-closed field, and let $F=k((X))$.  Now let $a_1=a_2=X+2$, so that $\alpha_1=\alpha_2=X(X+4)$. Note that $X+4$ is a non-zero square in $F$, since $4$ is a non-zero square.  In particular,  $(\alpha_1)=(\alpha_2)=(X)\in H^1(F)$.

	We now proceed to show that, for all $d\geq 2$, any normalized cohomological invariant $\eta\in \mathrm{Inv}_0^d(\textbf{GO}(b))$ which vanishes at $b_{-a_1}\otimes b_{-a_2}$ is identically zero.
	
	Let $d\geq 2$, and let $\eta\in   \mathrm{Inv}_0^d(\textbf{GO}(b))$.  By Theorem \ref{thmcoh}, we have $\eta=\rho\cup\nu$ for some $\rho\in H^{d-2}(F)$.
	Since $k$ is real-closed, one may write $$\rho=\varepsilon_0 (-1)^{\cup (d-2)}+\varepsilon_1 (-1)^{\cup (d-3)}\cup (X),$$ for some unique $\varepsilon_0,\varepsilon_1\in\mathbb{F}_2$.
	Hence, for any field extension $E/F$, and any $\lambda\in E^\times$, we have \[ \eta_E([ \lambda])=\left(\varepsilon_0 
	(-1)^{\cup (d-2)}+\varepsilon_1 (-1)^{\cup (d-3)}\cup (X)\right)\cup(\alpha_1)\cup(\lambda)\in H^d(E).\]
	Since $(X)\cup (\alpha_1)=(X)\cup(X)=(-1)\cup(X)$, one may assume without loss of generality that $\varepsilon_0=0$, so that  
	$\eta=\varepsilon (-1)^{(d-2)}\cup\nu$ for a unique $\varepsilon\in\mathbb{F}_2$.
	
	Now, we have \[\nu_F(b_{-a_1}\otimes b_{-a_2})=(\alpha_1,\alpha_2)=(-1)\cup(X).\] Thus, if $\eta$ vanishes at $b_{-a_1}\otimes b_{-a_2}$, then $\varepsilon (-1)^{d-1}\cup(X)=0$, and thus $\varepsilon=0$, implying in turn that $\eta=0$.

    To finish this example, 
    let us note that, for all $d\geq 2$, the invariant $\eta=(-1)^{d-2}\cup \nu$ satisfies $\eta_F\neq 0$, since the image of $b_{-a_1}\otimes b_{-a_2}$ under $\eta_F$ is $(-1)^{\cup(d-1)}\cup(X)$, which non-zero.
 \end{ex}

\begin{ex}\label{laurent}
 Let $k$ be an algebraically closed field of characteristic different from two, and let $F=\displaystyle\bigcup_{\substack{n\geq 1}}k((X_1))\cdots ((X_n))$.
For any subset $I=\{i_1,\ldots, i_n\}$ of $\mathbb{N}^*$ of cardinality $n$, set \[ (X_I)=(X_{i_1})\cup\cdots \cup (X_{i_n})\in H^n(F). \]

If $I=\emptyset$, we will use the convention $(X_\emptyset)=\overline{1}\in H^0(F)=\mathbb{F}_2$.

Note that the definition of $(X_I)$ does not depend on the choice of the numbering of the elements of $I$, since $H^*(F)$ is a commutative ring.

 Then, for all $n\geq 0$, $H^n(F)$ is an $\mathbb{F}_2$-vector space with basis \[ (X_I), I\subset\mathbb{N}^*, \vert I\vert =n. \]

 Now let $a_1=X_1+2, a_2=X_2+2$, so that $\alpha_i=X_i(X_i+4)$ for $i=1,2$. Note that $X_i+4$ is a non-zero square in $F$, since $4$ is a non-zero square, and thus $(\alpha_i)=(X_i)\in H^1(F)$ for $i=1,2$.

 We now proceed to show that, for all $d\geq 2$, any normalized cohomological invariant $\eta\in \mathrm{Inv}_0^d(\textbf{GO}(b))$ which vanishes at $b_{-a_1}\otimes b_{-a_2}$ is identically zero.

Let $\eta\in   \mathrm{Inv}_0^d(\textbf{GO}(b))$. Since $-1$ is a square in $k$, hence in $F$, every element of $F$ is a sum of two squares. By Theorem \ref{thmcoh}, we have $\eta=\rho\cup\nu$ for some $\rho\in H^{d-2}(F)$.

One may write  $\rho=\displaystyle\sum_{\vert I\vert=d-2}\varepsilon_I\cdot (X_I), \ \varepsilon_I\in\mathbb{F}_2$.

Note that, for any field extension $E/F$, and any  $u\in K_E^\times$, we have \[  (\alpha_1)\cup( N_{K_E/E}(u))=(\alpha_2 (\alpha_1\alpha_2))\cup ( N_{K_E/E}(u))=(\alpha_2)\cup( N_{K_E/E}(u))\in H^2(E), \] that is \[ (X_1)\cup( N_{K_E/E}(u))=(X_2)\cup( N_{K_E/E}(u))\in H^2(E). \]

Now, we have $(X_i)\cup(X_i)=(X_i)\cup(-1)=0\in H^2(E)$ since $-1$ is a square in $k$, hence in $E$. Thus, we get $(X_I)\cup \nu=0$  whenever $I$ contains $1$ or $2$.

Therefore, without any loss of generality, one may assume that $\varepsilon_I=0$ for all $I$ intersecting $\{1,2\}$ non-trivially, that is  \[ \rho=\sum_{\substack{I\subset \mathbb{N}_{\geq 3}\\ \vert I\vert =d-2}}\varepsilon_I\cdot (X_I). \]

Assume now  that $\eta_F(b_{-a_1}\otimes b_{-a_2})=0$. This means that \[ \sum_{\substack{I\subset \mathbb{N}_{\geq 3}\\ \vert I\vert =d-2}}\varepsilon_I\cdot (X_I)\cup(X_1)\cup(X_2)=0. \]

Since each $I$ in the sum is disjoint from $\{1,2\}$, we get \[ \sum_{\substack{I\subset \mathbb{N}_{\geq 3}\\ \vert I\vert =d-2}}\varepsilon_I\cdot (X_{I\cup\{1,2\}})=0, \] and thus $\varepsilon_I=0$ for all $I\subset \mathbb{N}_{\geq 3}$. Hence $\rho=0$, and therefore $\eta=0$.

Finally, note  that, for $d\geq 2$, the invariant $\eta=(X_{\{3,\ldots,d\}})\cup \nu\in \mathrm{Inv}_0^d(\textbf{GO}(b)) $ satisfies $\eta_F\neq 0$, since the image of $b_{-a_1}\otimes b_{-a_2}$ under $\eta_F$ is $(X_{\{1,\ldots,d\}})\neq 0\in H^d(F)$.
 \end{ex}

\end{document}